\documentclass[11pt]{amsart}
\usepackage{amsfonts,latexsym,amsthm,amssymb,graphicx}
\usepackage[all]{xy}

\DeclareFontFamily{OT1}{rsfs}{}
\DeclareFontShape{OT1}{rsfs}{n}{it}{<-> rsfs10}{}
\DeclareMathAlphabet{\mathscr}{OT1}{rsfs}{n}{it}

\CompileMatrices

\newtheorem{theorem}{Theorem}[section]
\newtheorem{lemma}[theorem]{Lemma}
\newtheorem{corol}[theorem]{Corollary}
\newtheorem{prop}[theorem]{Proposition}

{\theoremstyle{definition} \newtheorem{defin}[theorem]{Definition}}
{\theoremstyle{remark} \newtheorem{remark}[theorem]{Remark}
\newtheorem{example}{Example}[section]}

\newcommand{\Cbb}{{\mathbb{C}}}
\newcommand{\Qbb}{{\mathbb{Q}}}
\newcommand{\Gbb}{{\mathbb{G}}}

\newcommand{\Pbb}{{\mathbb{P}}}

\newcommand{\TPbb}{{\Til \Pbb}}

\newcommand{\PGL}{\text{\rm PGL}}
\newcommand{\Sym}{\text{\rm Sym}}

\newcommand{\cC}{{\mathscr C}}
\newcommand{\cB}{{\mathscr B}}
\newcommand{\cCred}{{\cC'}}
\newcommand{\OC}{{\mathscr O_\cC}}
\newcommand{\OCbar}{\overline{\OC}}

\newcommand{\cS}{{\mathscr S}}
\newcommand{\cX}{{\mathscr X}}

\newcommand{\Til}[1]{{\widetilde{#1}}}

\DeclareMathOperator{\rk}{rk}

\DeclareMathOperator{\im}{im}

\DeclareMathOperator{\Hom}{Hom}

\title[Limits of PGL(3)-translates of plane curves, I]
{Limits of PGL(3)-translates of plane curves, I}
\author{Paolo Aluffi, Carel Faber}
\address{Dept.~of Mathematics, Florida State University, Tallahassee
FL 32306, U.S.A.}
\email{aluffi@math.fsu.edu}
\address{Inst.~f\"or Matematik, Kungliga Tekniska H\"ogskolan, 
S-100 44 Stockholm, Sweden}
\email{faber@math.kth.se}

\begin{document}

\begin{abstract}
We classify all possible limits of families of translates of a fixed,
arbitrary complex plane curve. We do this by giving a set-theoretic 
description of the projective normal cone (PNC)
of the base scheme of a natural rational map, determined by
the curve, from the $\Pbb^8$ of $3\times 3$ matrices to the
$\Pbb^N$ of plane curves of degree~$d$.
In a sequel to this paper we determine the multiplicities of the 
components of the PNC. The knowledge of the PNC as a cycle is 
essential in our computation of the
degree of the $\PGL(3)$-orbit closure of an arbitrary plane curve, 
performed in~\cite{MR2001h:14068}. 
\end{abstract}

\maketitle

\section{Introduction}\label{intro}

In this paper we determine the possible {\em limits\/} of a fixed, 
arbitrary complex plane curve $\cC$, obtained by applying to it a family 
of translations $\alpha(t)$ centered at a singular transformation of the 
plane. In other words, we describe the curves in the boundary of the
$\PGL(3)$-orbit closure of a given curve $\cC$.

Our main motivation for this work comes 
from enumerative geometry. In 
\cite{MR2001h:14068}
we have determined the {\em degree\/} of the $\PGL(3)$-orbit closure 
of an arbitrary (possibly singular, reducible, non-reduced) plane curve;
this includes as special cases the determination of
several characteristic numbers of families of plane curves, 
the degrees of certain maps to moduli spaces of plane curves, and 
isotrivial versions of the Gromov-Witten invariants of the plane. 
A description of the limits of a curve, and in fact a more refined
type of information is an 
essential ingredient of our approach. This information
is obtained in this paper and in its sequel \cite{ghizzII}; the results
were announced and used in~\cite{MR2001h:14068}.

The set-up is as follows. Consider the natural action of $\PGL(3)$
on the projective space of plane curves of a fixed degree.
The orbit closure of a 
curve $\cC$ is dominated by the closure~$\TPbb^8$ of the
graph of the rational map $c$ from the $\Pbb^8$ of $3\times3$ 
matrices to the $\Pbb^N$ of plane curves of degree~$d$,
associating to $\varphi\in \PGL(3)$ the translate of $\cC$ by~$\varphi$. 
The boundary of the orbit consists of
limits of $\cC$ and plays an important role in the study of the orbit
closure.

Our computation of the degree of the orbit closure of $\cC$ 
hinges on the study of~$\TPbb^8$, and
especially of the scheme-theoretic inverse image in $\TPbb^8$ of the 
base scheme~$\cS$ of~$c$. 
Viewing $\TPbb^8$ as the blow-up of $\Pbb^8$ along $\cS$, this
inverse image is the exceptional divisor, and may be identified
with the projective normal cone (PNC) of $\cS$ in $\Pbb^8$.
A description of the PNC leads to a description of the limits 
of $\cC$: the image of the PNC in $\Pbb^N$ is contained in the set of 
limits, and the complement, if nonempty, consists of easily identified 
`stars' (that is, unions of concurrent lines).

This paper is devoted to a set-theoretic description of the PNC for an
arbitrary curve. This
suffices for the determination of the limits, but does not
suffice for the enumerative applications in~\cite{MR2001h:14068};
these applications require the full knowledge of the PNC {\em as a cycle,\/} 
that is, the determination of the multiplicities of its different components. 
We obtain this additional information in \cite{ghizzII}. 

The final result of our analysis (including multiplicities)
was announced in \S2 of~\cite{MR2001h:14068}.
The proofs of the facts stated there are given in the present article and its
sequel.
The main theorem of this paper (Theorem~\ref{mainmain}, in \S\ref{proof})
gives a precise set-theoretic description of the PNC, relying upon five
types of families and limits identified in \S\ref{germlist}. 
In this introduction 
we confine ourselves to 
formulating a weaker version, focusing on the determination of limits.
In~\cite{ghizzII} (Theorem~2.1), we compute the multiplicities of the 
corresponding five types of components of the PNC.

The limits of a curve $\cC$ are necessarily curves with {\em small 
linear orbit,\/} that is, curves with infinite stabilizer.
Such curves are classified in \S1 of~\cite{MR2002d:14084};
we reproduce the list of curves obtained in \cite{MR2002d:14084}
in an appendix at the end of this paper (\S\ref{appendix}).
For another classification, from a somewhat different viewpoint,
we refer to \cite{MR1698902}.
For these curves, the limits can be determined 
using the results in~\cite{MR2002d:14083} 
(see also~\S\ref{boundary}). The following statement reduces the computation 
of the limits of an arbitrary curve $\cC$ to the case of curves with small orbit.

\begin{theorem}\label{main}
Let $\cX$ be a limit of a plane curve $\cC$ of degree~$d$, obtained by 
applying to it a $\Cbb((t))$-valued point of $\PGL(3)$ with singular center.
Then $\cX$ is in the orbit closure of a
star (reproducing projectively the $d$-tuple cut out on 
$\cC$ by a line meeting it properly), 
or of curves with small orbit determined by
the following features of~$\cC$:
\begin{itemize}
\item[I:] The linear components of the support $\cCred$ of $\cC$;
\item[II:] The nonlinear components of $\cCred$;
\item[III:] The points at which the tangent cone of $\cC$ is 
supported on at least $3$ lines;
\item[IV:] The Newton polygons of $\cC$ at the singularities and 
inflection points of $\cCred$;
\item[V:] The Puiseux expansions of formal branches of $\cC$ at the 
singularities of $\cCred$.
\end{itemize}
\end{theorem}

The limits corresponding to these features may be described as follows.
In cases I and III they are 
unions of a star and a general line, that we call `fans';
in case II, they are supported on the union of a nonsingular conic and
a tangent line; 
in case IV, they are supported on the union of the coordinate triangle
and several curves from a pencil $y^c=\rho\, x^{c-b} z^b$,
with $b<c$ coprime positive integers; and in case~V
they are supported on unions of quadritangent conics and 
the distinguished tangent line.
The following picture illustrates the limits in cases~IV and~V:
\begin{center}
\includegraphics[scale=.6]{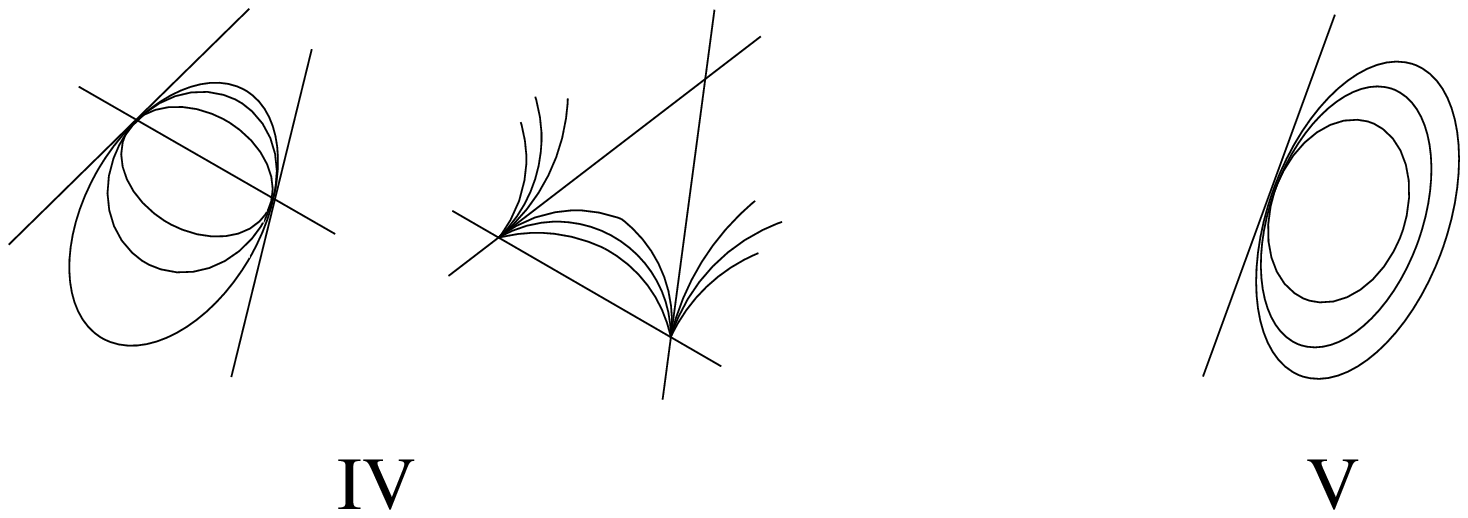}
\end{center}
A more precise description of the limits is given in \S\ref{germlist}, 
referring to the classification of these curves obtained in \S1 of 
\cite{MR2002d:14084}
and reproduced in \S\ref{appendix} of this paper.

The proof of Theorem~\ref{main} (or rather of its more precise form
given in Theorem~\ref{mainmain}) is by an explicit reduction process, 
and goes along the following lines. 
The stars mentioned in the 
statement are obtained by families of translations $\alpha(t)$
(`germs') centered at an element $\alpha(0)\not\in\cS$.
To analyze germs centered at points of $\cS$, 
we introduce a notion of equivalence of germs
(Definition~\ref{equivgermsnew}), such that equivalent germs 
lead to the same limit. We then prove that every germ centered
at a point of~$\cS$ is essentially
equivalent to one with matrix representation
$$\begin{pmatrix}
1 & 0 & 0\\
q(t) & t^b & 0\\
r(t) & s(t)t^b & t^c
\end{pmatrix}$$
with $0\le b\le c$ and $q$, $r$, and $s$ polynomials. Here,
coordinates are chosen so that the point $p=(1:0:0)$ belongs to $\cC$. 
Studying the limits obtained by applying such germs to $\cC$,
we identify five specific types of families 
(the {\em marker germs\/} listed in~\S\ref{germlist}), reflecting the 
features of $\cC$ at $p$ listed in Theorem~\ref{main}, and with the 
stated kind of limit. We prove that unless the germ is of one of
these types, the corresponding limit is already accounted for (for 
example, it is in the orbit closure of a star of the type mentioned in the statement).

In terms of the graph of the rational map $c$ mentioned above, we
prove that every component of the PNC is hit at a general point 
by the lift in $\TPbb^8$ of 
one of the five distinguished types of germs. This yields our set-theoretic 
description of the PNC. In fact, the lifts intersect the 
corresponding components transversally, and this will be important 
in our determination of the multiplicities of the components in~\cite{ghizzII}.

The procedure underlying the proof of Theorem~\ref{mainmain} may
be applied to any given plane curve, producing a list of its limits.
In practice, one needs to find the marker germs for the curve;
these determine the components of the PNC. The two examples 
in \S\ref{twoexamples} illustrate this process, and show that components
of all types may already occur on curves of degree~$4$.
Here is a simpler example, for a curve of degree~$3$.

\begin{example}\label{exone}
Consider the irreducible cubic~$\cC$ given by the equation
$$xyz+y^3+z^3=0\,.$$
It has a node at $(1:0:0)$ and three inflection points.
According to Theorem~\ref{mainmain} and the list in \S\ref{germlist},
the PNC for~$\cC$ has one component of type~II and several of 
type~IV. The latter correspond to the three inflection points and
the node. A list of representative marker germs for the component of
type~II and for the component of type~IV due to the node 
may be obtained by following the procedure explained in \S\ref{setth}:
\[
{\rm II}:
\begin{pmatrix}
-2 & -t & 0\\
1 & t & 0\\
1 & 0 & t^2
\end{pmatrix};\quad
{\rm IV}:
\begin{pmatrix}
1 & 0 & 0\\
0 & t & 0\\
0 & 0 & t^2
\end{pmatrix}
\,,\,
\begin{pmatrix}
1 & 0 & 0\\
0 & t^2 & 0\\
0 & 0 & t
\end{pmatrix}\,.
\]
The latter two marker germs, corresponding to the two lines in the tangent
cone at the node, have the same center and lead to 
projectively equivalent limits, hence they contribute the same component 
of the PNC. 
Equations for the limits of $\cC$ determined by the germs listed above
are
$$
x(xz+2y^2)=0,\quad
y(y^2+xz)=0,\quad \text{and} \quad
z(z^2+xy)=0\,,
$$
respectively: a conic with a tangent line,
and a conic with a transversal line
(two limits). The inflection points also contribute components
of type~IV; the limits in that case are cuspidal cubics.

According to Theorem~\ref{main}, all limits of $\cC$ (other than stars of
lines) are projectively equivalent to one of these curves, or to limits of them
(cf.~\S\ref{boundary}).
\end{example}

Necessary preliminary considerations, and the full statement of the 
main theorem, are found in \S\ref{prelim}.
The determination of the limits by successive reductions of a given 
family of curves, proving the result, is worked out in \S\ref{setth}
and~\S\ref{typeVcomps}. In
\S\ref{boundary} we summarize the more straightforward situation 
for curves with small orbits.\smallskip

Harris and Morrison (\cite{MR99g:14031}, p.~138) pose the {\em flat 
completion problem\/} for families of embedded curves, asking for
the determination of all curves in $\Pbb^n$ that can arise as flat limits
of a family of embedded stable curves over the punctured disc.
The present article solves the isotrivial form of this problem, for 
plane curves.

In principle, a solution of the isotrivial flat completion problem for plane
curves can already be found in the marvelous article~\cite{ghizz}
by Aldo Ghizzetti, dating back to the 1930s. However,
Ghizzetti's
results do not lead to a description of the PNC, which is 
necessary for our application in \cite{MR2001h:14068}, and which 
is the main result of this paper and of its sequel.

Caporaso and Sernesi use our determination of the limits in
\cite{MR2003k:14035} (Theorem 5.2.1). Hacking \cite{MR2078368}
and Hassett \cite{MR2000j:14045} study the limits of families of
nonsingular plane curves of a given degree, by methods different from
ours: they allow the plane to degenerate together with the curve. It
would be interesting to compare their results to ours.
However, there are fundamental differences between the phenomena we
study and those addressed in \cite{MR2078368} and \cite{MR2000j:14045}; 
for example, our families are constant in moduli, 
and our results apply to {\em arbitrary\/} plane curves.
By the same token, neither Hacking-stability nor GIT-stability
play an important role in our study. 
Consider the case of a plane curve with an analytically 
irreducible singularity. The determination of the contribution of the
singularity to the PNC of the curve requires both its {\em linear\/}
type and {\em all\/} its Puiseux pairs, see \S5 of \cite{MR2001h:14068}.
In general, the stability
conditions mentioned above require strictly less 
(cf.~Kim-Lee \cite{MR2090618}). 
For example, a singularity analytically isomorphic to $y^2=x^5$
on a {\em quartic\/} leads necessarily to a component of type~V
(cf.~Example~\ref{extwo}), whereas on a quintic, it leads to
either a component of type~IV or a component of type~V,
according to the order of contact with the tangent line.
For GIT-stability, see also Remark~\ref{GITrem}.

\smallskip
The enumerative problem considered in \cite{MR2001h:14068}, as 
well as the question of limits of PGL-translates, makes sense for
hypersurfaces of projective space of any dimension. The case of 
configurations of points in $\Pbb^1$ is treated in \cite{MR1244973}.
The degree of the orbit closure of a configuration of planes in $\Pbb^3$
is computed in \cite{MR2455792}. In general, these problems appear to be
very difficult. The techniques used in this paper could in principle be
used in arbitrary dimension, but the case-by-case analysis (which is 
already challenging for curves in $\Pbb^2$) would likely be unmanageable
in higher dimension. By contrast, the techniques developed in \cite{ghizzII}
should be directly applicable: once `marker germs' have been determined,
computing the multiplicities of the corresponding components of the PNC
should be straightforward,
using the techniques of \cite{ghizzII}.

\vskip6pt

\thanks{{\bf Acknowledgments.} Work on this paper was made
possible by support from 
Mathematisches Forschungsinstitut Oberwolfach,
the Volkswagen Stiftung,
the Max-Planck-Institut f\"ur Mathematik (Bonn),
Princeton University,
the G\"oran Gustafsson foundation,
the Swedish Research Council,
the Mittag-Leffler Institute,
MSRI, NSA, NSF, and our home institutions. 
We thank an anonymous 
referee of our first article on the topic of linear orbits of plane curves,
\cite{MR94e:14032}, for bringing the paper of Aldo Ghizzetti to our
attention.} 
We also thank the referee of this paper and \cite{ghizzII}, for the
careful reading of both papers and for comments
that led to their improvement.

\section{Set-theoretic description of the PNC}\label{prelim}

\subsection{Limits of translates}\label{rough}
We work over $\Cbb$. 
We choose homogeneous coordinates $(x:y:z)$ in $\Pbb^2$, and 
identify $\PGL(3)$ with the open set of nonsingular matrices in the 
space $\Pbb^8$ parametrizing $3\times 3$ matrices. 
We consider the right action of  $\PGL(3)$ on the space 
$\Pbb^N=\Pbb H^0(\Pbb^2, \mathcal O(d))$ of degree-$d$ plane 
curves; if $F(x,y,z)=0$ is an equation for a plane curve~$\cC$, 
and $\alpha\in \PGL(3)$, we denote by $\cC\circ\alpha$ the curve 
with equation $F(\alpha(x,y,z))=0$.

We will consider families of plane curves over the punctured disk, of
the form $\cC\circ\alpha(t)$, where $\alpha(t)$ is a $3\times 3$ matrix
with entries in $\Cbb[t]$, such that $\alpha(0)\ne 0$, $\det\alpha(t)
\not\equiv 0$, and $\det\alpha(0)=0$.
Simple reductions show that studying these families is equivalent to
studying all families $\cC\circ\alpha(t)$, where $\alpha(t)$ is a 
$\Cbb((t))$-valued point of~$\Pbb^8$ such that $\det\alpha(0)=0$.
We also note that if $\cC$ is a smooth curve of degree $d\ge 4$,
then any family of curves of degree~$d$ parametrized by the punctured
disk and whose members are abstractly isomorphic to $\cC$, i.e., an
isotrivial family,
is essentially of this type (cf.~\cite{MR770932}, p.~56).

The arcs of matrices $\alpha(t)$ will be called {\em germs,\/} and viewed 
as germs of curves in~$\Pbb^8$. 
The flat limit $\lim_{t\to 0}\,\cC\circ \alpha(t)$
of a family $\cC\circ\alpha(t)$ as $t \to 0$
may be computed concretely by clearing common powers of
$t$ in the expanded expression $F(\alpha(t))$, and then setting $t=0$.
Our goal is the determination of all possible limits of families as above,
for a given arbitrary plane curve~$\cC$.

\subsection{The Projective Normal Cone}\label{ident}
The set of all translates $\cC\circ\alpha$ is the {\em linear orbit\/} of
$\cC$, which we denote by $\OC$; the complement of $\OC$ in its
closure $\OCbar$ is the {\em boundary\/} of the orbit of $\cC$. 
By the {\em limits of $\cC$\/} we will mean the limits of families 
$\cC\circ\alpha(t)$ with $\alpha(0)\not\in \PGL(3)$.

\begin{remark}
For every curve $\cC$, the boundary is a subset of the set of limits; 
if $\dim\OC=8$ (the stabilizer of $\cC$ is finite), then these two sets
coincide. If $\dim\OC<8$ (the stabilizer is infinite, and the orbit is 
{\em small,\/} in the terminology of \cite{MR2002d:14083} and 
\cite{MR2002d:14084}) then there are families with limit equal 
to~$\cC$; in this case, the whole orbit closure $\OCbar$ consists of
limits of $\cC$.
\end{remark}

The set of limit curves is itself a union of orbits of plane curves; 
our goal is a description of representative elements of these orbits; 
in particular, this will yield a description of the boundary of $\OC$.
In this section we relate the 
set of limits of $\cC$ to the {\em projective normal cone\/} mentioned 
in the introduction.

Points of $\Pbb^8$, that is, $3\times 3$ matrices, may be viewed
as rational maps $\Pbb^2 \dashrightarrow \Pbb^2$. 
The kernel of a singular matrix $\alpha\in \Pbb^8$ determines a 
line of $\Pbb^2$ (if $\rk\alpha=1$) or a point
(if $\rk\alpha=2$); $\ker\alpha$
will denote this locus. Likewise, the image of $\alpha$ is a point of
$\Pbb^2$ if $\rk\alpha=1$, or a line if $\rk\alpha=2$.

The action map
$\alpha \mapsto \cC\circ\alpha$
for $\alpha\in \PGL(3)$ defines a rational map
$$c: \Pbb^8 \dashrightarrow \Pbb^N\quad.$$
We denote by $\cS$ the base scheme of this rational map.
The closure of the graph of $c$ may be identified with the blow-up 
$\TPbb^8$ of $\Pbb^8$ along $\cS$. 
The support of $\cS$ consists of the matrices $\alpha$ such that 
(with notation as above) $F(\alpha(x,y,z))\equiv 0$; that is, matrices 
whose image is contained in $\cC$.

The {\em projective normal cone\/} (PNC) of $\cS$ in $\Pbb^8$ is the
exceptional divisor $E$ of this blow-up.
We have the following commutative diagram:
$$\xymatrix@M=10pt{
E \ar@{->>}[d] \ar@<-.5ex>@{^(->}[r] & \TPbb^8 \ar@{->>}[d]^\pi 
\ar@<-.5ex>@{^(->}[r] & \Pbb^8\times \Pbb^N   \ar@{->>}[d] \\
\cS \ar@<-.5ex>@{^(->}[r] & \Pbb^8 \ar@{-->}[r]^c & \Pbb^N
}$$
Therefore, as a subset of $\Pbb^8\times\Pbb^N$, the support of 
the PNC is
\begin{multline*}
|E|=\{(\alpha,\cX)\in \Pbb^8\times\Pbb^N : \text{$\cX$ is a limit of
  $\cC\circ \alpha(t)$}\\
\text{for some germ $\alpha(t)$ centered at
  $\alpha\in \cS$ and not contained in $\cS$}\}\quad.
\end{multline*}

\begin{lemma}\label{PNCtolimits}
The set of limits of $\cC$ consists of the image of the PNC in $\Pbb^N$,
and of limits of families $\cC\circ \alpha(t)$ with $\alpha=\alpha(0)$ a 
singular matrix whose image is not contained in~$\cC$. 

In the latter case: if $\alpha$ has rank~1, the limit consists of a 
multiple line supported on $\ker\alpha$; if $\alpha$ has rank~2, 
the limit consists of a star of lines through $\ker\alpha$, reproducing 
projectively the tuple of points cut out by $\cC$ on the image of $\alpha$.
\end{lemma}

\begin{proof}
The PNC dominates the set of limits
of families $\cC\circ \alpha(t)$ for which $\alpha(t)$ is centered at
a point of indeterminacy of $c$. This gives the first statement.

To verify the second assertion, assume that $\alpha(t)$ is centered
at a singular matrix $\alpha$ at which $c$ {\em is\/} defined; $\alpha$ is 
then a rank-1 or rank-2 matrix such that $F(\alpha(x,y,z))\not\equiv 0$. 
After a coordinate change we may assume without loss of generality that
$$\alpha=\begin{pmatrix}
1 & 0 & 0 \\
0 & 0 & 0 \\
0 & 0 & 0
\end{pmatrix}
\quad\text{or}\quad
\alpha=\begin{pmatrix}
1 & 0 & 0 \\
0 & 1 & 0 \\
0 & 0 & 0
\end{pmatrix}$$
and $F(x,0,0)$, resp.~$F(x,y,0)$ are not identically zero. These are
then the forms defining the limits of the corresponding families, and the 
descriptions given in the statement are immediately verified in these 
cases.
\end{proof}

The second part of Lemma~\ref{PNCtolimits} may be viewed as the 
analogue in our context of an observation of Pinkham (`sweeping out 
the cone with hyperplane sections', \cite{MR0376672},~p.~46).

\begin{remark}\label{eluding}
Denote by $R$  the proper transform in $\TPbb^8$ of the set of singular
matrices in $\Pbb^8$. Lemma~\ref{PNCtolimits} asserts that
the set of limits of $\cC$ is the image of the union of the PNC and $R$.
A more explicit description of the image
of $R$ has eluded us; 
for a smooth curve~$\cC$ of degree $\ge 5$ these `star limits' have 
two moduli. It would be interesting to obtain a classification of curves 
$\cC$ with smaller `star-moduli'.

The image of the {\em intersection\/} of $R$ and the PNC will play an
important role in this paper. Curves in the image of this locus will
be called `rank-$2$ limits'; we note that the set of rank-$2$ limits
has dimension~$\le 6$.
\end{remark}

Lemma~\ref{PNCtolimits} translates the problem of finding the limits
for families of plane curves $\cC\circ\alpha(t)$ into the problem 
of describing the PNC for the curve $\cC$.
Each component of the PNC is a $7$-dimensional irreducible subvariety 
of $\TPbb^8\subset \Pbb^8\times \Pbb^N$. We will describe it by listing 
representative points of the component. 
More precisely, note that $\PGL(3)$ acts on $\Pbb^8$ by right multiplication, 
and that this action lifts to a right action of $\PGL(3)$ on $\TPbb^8$. 
Each component of the PNC is a union of orbits of this action.
For each component, we will list germs $\alpha(t)$ lifting 
on~$\TPbb^8$ to germs $\tilde\alpha(t)$ so that the union of the orbits of 
the centers $\tilde\alpha(0)$ is dense in that component.

\subsection{Marker germs}\label{germlist}
In a coarse sense, the classification of limits into `types' as in 
Theorem~\ref{main} depends on the image of the center $\alpha(0)$ of 
the family: this will be a subset of $\cC$ (cf.~Lemma~\ref{PNCtolimits}), 
hence it will either be a (linear) component of $\cC$ (type~I), or a point 
of $\cC$ (general for type~II, singular or inflectional for types III, IV, and~V).

We will now list germs determining the components of the PNC in the
sense explained above. We will call such a germ a {\em marker germ,\/}
as the center of its lift to~$\TPbb^8$ 
(the corresponding {\em marker center\/}) `marks' a component of the~PNC.
The first two types depend on global features of $\cC$:
its linear and nonlinear components. The latter three
depend on local features of $\cC$: inflection points and
singularities of (the support of) $\cC$. 
That there are only two global types is due to the fact that
the order of contact of a nonlinear component and the
tangent line at a general point equals two (in characteristic zero).
The three local types are due to linear features at
singularities of $\cC$ (type~III), single nonlinear branches at 
special points of $\cC$ (type~IV), and collections of several 
matching nonlinear branches at singularities of $\cC$ (type~V).
Only type~V leads to limits with additive stabilizers, and the absence
of further types is due to the fact, shown in \cite{MR2002d:14084}, 
that in characteristic zero only one kind of curves with small orbit 
has additive stabilizers (also cf.~\S\ref{appendix}).

\begin{remark}\label{GITrem}
A plane curve with small orbit is not GIT-stable. Whether it is strictly
semistable or unstable is not directly related to the questions we are 
considering here. For example, the curves $xyz$ and $x^2yz$ have 
similar behavior from the point of view of this paper; yet the former is
strictly semistable, the latter is unstable.

Similarly, consider the union of a general quartic and a multiple
line in general position. This has 8-dimensional orbit; it is stable
in degree~5, strictly semistable in degree~6, and unstable in higher
degrees. But the multiplicity of the line does not affect the behavior
from our point of view in any substantial way.

The lesson we draw from these examples is that there is no direct
relation between the considerations in this paper and GIT. We should
point out that the referee of this paper suggests otherwise, noting that
closures of orbits are of interest in both contexts, curves with small
orbits play a key role, and the mechanics of finding the limits is
somewhat similar in the two situations. The referee asks: {\em which
marker germs would be relevant in a GIT analysis?\/} We pass this
question on to the interested reader.
\end{remark}

The terminology
employed in the following matches the one in \S2 of \cite{MR2001h:14068};
for example, a {\em fan\/} is the union of a star and a general line.
In four of the five types, $\alpha=\alpha(0)$ is a rank-1 matrix and the 
line $\ker\alpha$ plays an important role; we will call this `the kernel
line'.

{\bf Type I.\/}
Assume $\cC$ contains a line, defined by a linear polynomial $L$. 
Write a generator of the ideal of $\cC$ as
$$F(x,y,z)=L(x,y,z)^m G(x,y,z)$$
with $L$ not a factor of $G$. Type I limits are obtained by germs
$$\alpha(t)=\alpha(0)+t\beta(t)\quad,$$
where $\alpha(0)$ has rank~2 and image the line defined by $L$. 

As we are assuming (cf.~\S\ref{rough}) that 
$\det\alpha(t)\not\equiv 0$, the
image of  $\beta(t)$ is not contained in $\im\alpha(0)$, so that the limit
$\lim_{t\to 0}L\circ\beta(t)$ is a well-defined line $\ell$. 
The limit $\lim_{t\to 0} \cC\circ\alpha(t)$ consists
of the $m$-fold line $\ell$, and a star of lines through the point
$\ker\alpha(0)$. This star reproduces projectively the tuple cut out on
$L$ by the curve defined by $G$.
\begin{center}
\includegraphics[scale=.4]{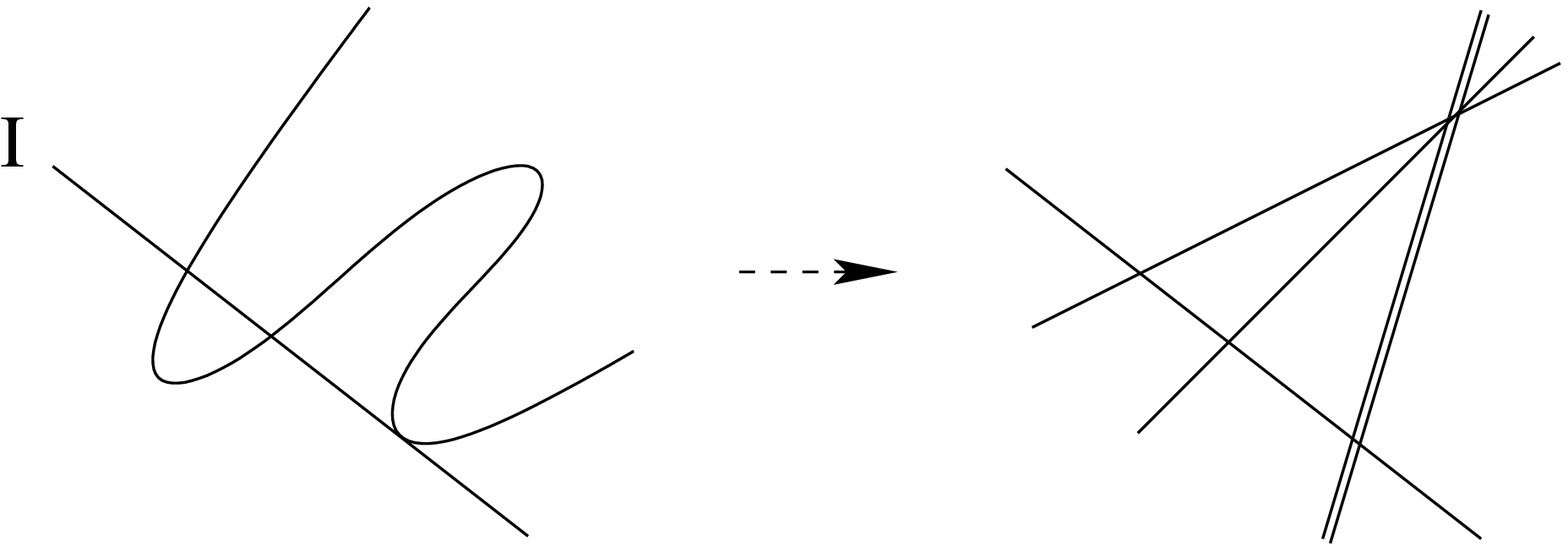}
\end{center}

The limit is in general a fan, and degenerates to a star if the $m$-fold 
line $\ell$ contains the point $\ker\alpha(0)$.
Fans and stars are studied in \cite{MR2002d:14084},
and are the only kinds of curves with small orbit that consist of lines; 
they are items (1) through (5) in our classification of curves with small 
orbit, see \S\ref{appendix}.\smallskip

For types~II---V we choose coordinates so that $p=(1:0:0)$ is a point
of $\cC$; for types II, IV, and V we further require that $z=0$ is a
chosen component $\ell$ of the tangent cone to $\cC$ at $p$.

{\bf Type II.\/}
Assume that $p$ is a nonsingular, non-inflectional point of 
the support $\cCred$ of $\cC$, contained in a nonlinear component,
with tangent line $z=0$. 
Let
$$\alpha(t)=\begin{pmatrix}
1 & 0 & 0 \\
0 & t & 0 \\
0 & 0 & t^2
\end{pmatrix}\quad.$$
Then the ideal of $\lim_{t\to 0}\cC\circ \alpha(t)$ is generated by
$$x^{d-2S}(y^2+\rho x z)^S\quad,$$ 
where $S$ is the multiplicity of the component in $\cC$, and $\rho\ne 0$;
that is, the limit consists of a 
(possibly multiple) nonsingular conic tangent to the kernel line,
union (possibly) a multiple of the kernel line.

\begin{center}
\includegraphics[scale=.5]{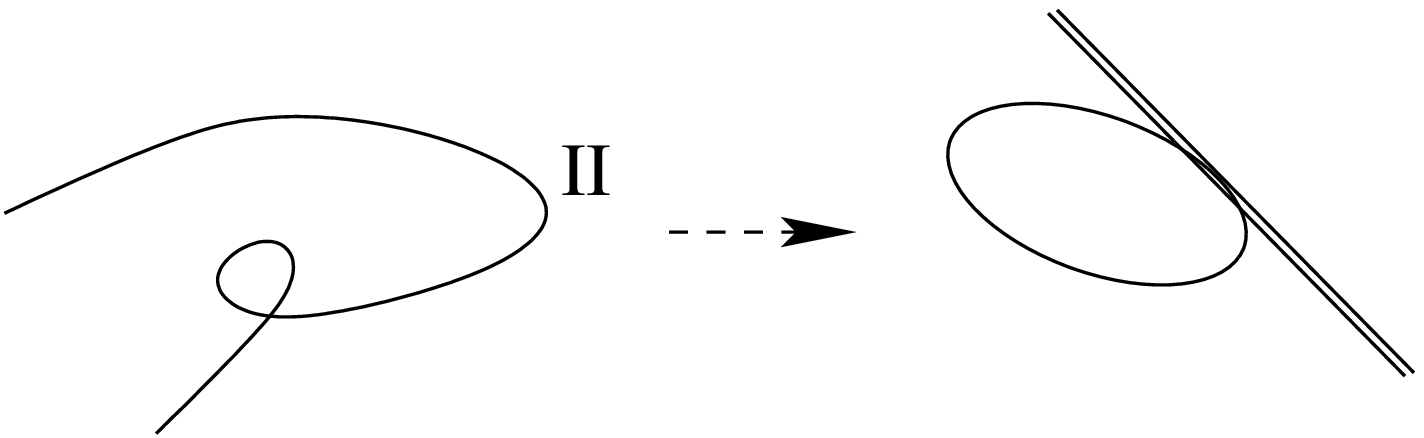}
\end{center}

Such curves are items (6) and (7) in the classification reproduced in 
\S\ref{appendix}.
The extra kernel line is present precisely when $\cC$ is
not itself a multiple nonsingular conic.

{\bf Type III.\/}
Assume that $p$ is a singular point of $\cCred$
of multiplicity $m$ in $\cC$, with tangent cone supported on at least three lines. 
Let
$$\alpha(t)=\begin{pmatrix}
1 & 0 & 0 \\
0 & t & 0 \\
0 & 0 & t
\end{pmatrix}\quad.$$
Then $\lim_{t\to 0} \cC\circ\alpha(t)$ is a fan consisting of a star 
centered at $(1:0:0)$ and projectively equivalent to the tangent cone to 
$\cC$ at $p$, and of a residual $(d-m)$-fold line supported on the 
kernel line $x=0$.
\begin{center}
\includegraphics[scale=.5]{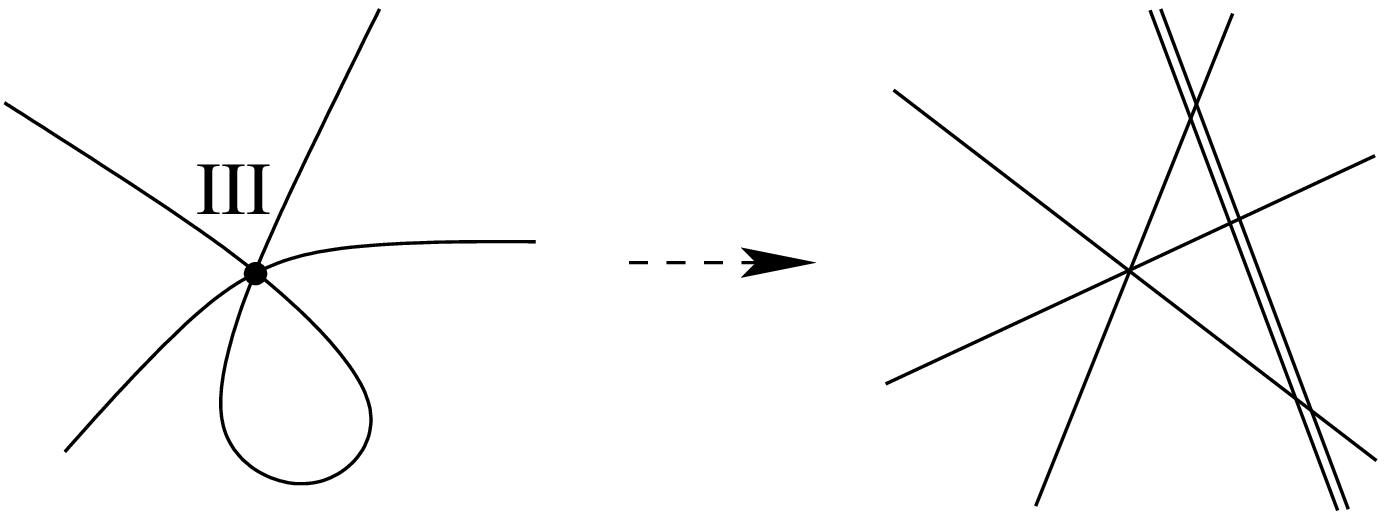}
\end{center}

{\bf Type IV.\/}
Assume that $p$ is a singular or inflection point of the support 
of $\cC$.
Germs of type~IV are determined by the choice of the line $\ell$ in the tangent 
cone to $\cC$ at $p$, and by the choice of a side of a corresponding 
Newton polygon, with slope strictly between $-1$ and $0$. 
This procedure is explained in more detail in \S\ref{details}.

Let $b<c$ be relatively prime positive integers such that  $-b/c$ is 
the slope of the chosen side. Let 
$$\alpha(t)=\begin{pmatrix}
1 & 0 & 0 \\
0 & t^b & 0 \\
0 & 0 & t^c
\end{pmatrix}\quad.$$
Then the ideal of $\lim_{t\to 0} \cC\circ\alpha(t)$ is generated by a
polynomial of the form
$$x^{\overline e}y^fz^e \prod_{j=1}^S(y^c+\rho_j x^{c-b}z^b)\quad,$$
with $\rho_j \ne 0$.
The number $S$ of `cuspidal' factors in the limit curve
is the number of segments cut out by the integer lattice on the
selected side of the Newton polygon.
\begin{center}
\includegraphics[scale=.5]{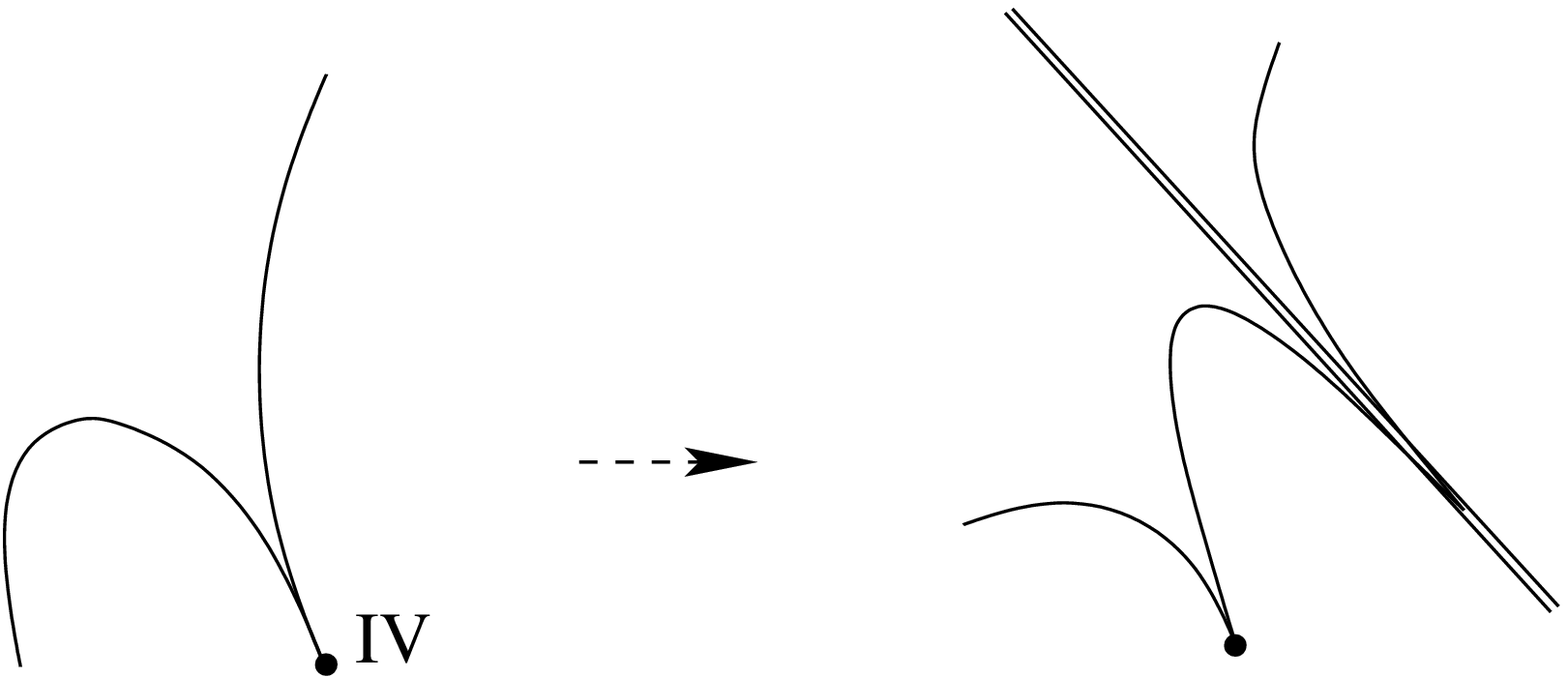}
\end{center}

The germ listed above contributes a component of the PNC unless
$b/c=1/2$ and the limit curve is supported on a conic union (possibly)
the kernel line.
The limit curves arising in this way are items (7) through (11)
listed in \S\ref{appendix}. (In particular, the picture drawn above does
not capture the possible complexity of the situation: several cuspidal
curves may appear in the limit, as well as all lines of the basic triangle.)
These limit curves are studied enumeratively in \cite{MR2002d:14083}.
The limit curves contributing components to the PNC in this fashion are 
precisely the curves that contain nonlinear components 
and for which the maximal connected subgroup of the stabilizer 
of the union of the curve and the kernel line
is the multiplicative group $\Gbb_m$.

{\bf Type V.\/}
Assume $p$ is a singular point of the support of $\cC$.
Germs of type~V are determined by the choice of the line $\ell$ in the 
tangent cone to $\cC$ at $p$, the choice of a formal branch
$z=f(y)=\gamma_{\lambda_0}y^{\lambda_0}+\dots$ for $\cC$ at $p$ 
tangent to~$\ell$, and the choice of a certain `characteristic' rational 
number $C>\lambda_0$ (assuming these choices can be made). 
This procedure is also explained in more detail in \S\ref{details}.

For $a<b<c$ positive integers such that 
$\frac ca=C$ and $\frac ba=
\frac{C-\lambda_0}2+1$, let
$$\alpha(t)=\begin{pmatrix}
1 & 0 & 0 \\
t^a & t^b & 0 \\
\underline{f(t^a)} & \underline{f'(t^a)t^b} & t^c
\end{pmatrix}$$
where $\underline{\cdots}$ denotes the truncation modulo $t^c$. 
The integer $a$ is chosen to be the minimal one
for which all entries in this germ are polynomials.
Then
$\lim_{t\to 0}\cC\circ \alpha(t)$ is given by
$$x^{d-2S}\prod_{i=1}^S\left(zx-\frac {\lambda_0(\lambda_0-1)}2
\gamma_{\lambda_0}y^2 -\frac{\lambda_0+C}2
\gamma_{\frac{\lambda_0+C}2}yx-\gamma_C^{(i)}x^2\right)\quad,$$
where $S$ and $\gamma_C^{(i)}$ are defined in \S\ref{details}.
\begin{center}
\includegraphics[scale=.4]{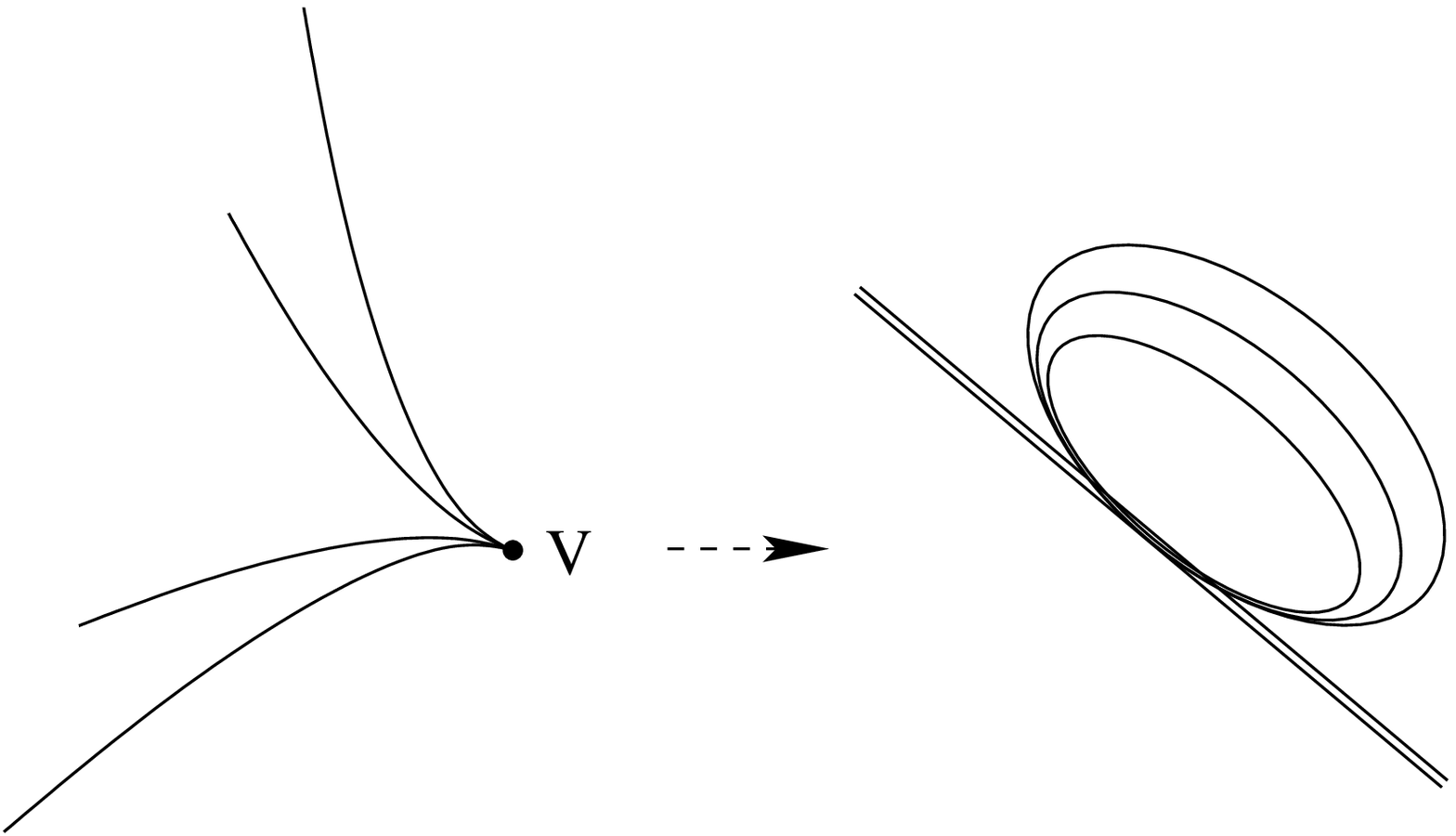}
\end{center}

These curves consist of {\em at least two\/} `quadritangent' conics---that is, 
nonsingular conics meeting at exactly one point---and (possibly)
a multiple kernel line.
(Again, the picture drawn here does not capture the subtlety of
the situation: these limits may occur already for irreducible
singularities.)
These curves are item (12) in the list in \S\ref{appendix},
and are studied enumeratively in \cite{MR2002d:14083},~\S4.1.
They are precisely the curves for which the maximal
connected subgroup of the stabilizer is the additive group $\Gbb_a$.

\subsection{Details for types IV and V}\label{details}
{\em Type IV:\/} Let $p=\im\alpha(0)$ be a singular or inflection point
of the support of $\cC$; choose a line in the tangent cone to
$\cC$ at $p$, and choose coordinates $(x:y:z)$ 
as before, so that $p=(1:0:0)$ and the selected line in the
tangent cone has equation $z=0$.
 The {\em Newton polygon\/} for $\cC$ in the
chosen coordinates is the boundary of the convex hull of the union of
the positive quadrants with origin at the points $(j,k)$ for which the
coefficient of $x^iy^jz^k$ in the generator~$F$ for the ideal of
$\cC$ is nonzero (see
\cite{MR88a:14001}, p.~380). The part of the Newton polygon consisting
of line segments with slope strictly between $-1$ and~$0$ does not
depend on the choice of coordinates fixing the flag $z=0$, $p=(1:0:0)$.

The limit curves are then obtained by choosing a side of the 
polygon with slope strictly between $-1$ and $0$, and setting to~$0$ 
the coefficients of the monomials in $F$ {\em not\/} on that side. 
These curves are
studied in \cite{MR2002d:14083}; typically, they consist of a union of
cuspidal curves. The kernel line is part of the distinguished triangle
of such a curve, and in fact it must be one of the distinguished tangents.

Here is the Newton polygon for the curve of Example~\ref{exone}, 
with respect to the point $(1:0:0)$ and the line $z=0$:
\begin{center}
\includegraphics[scale=.4]{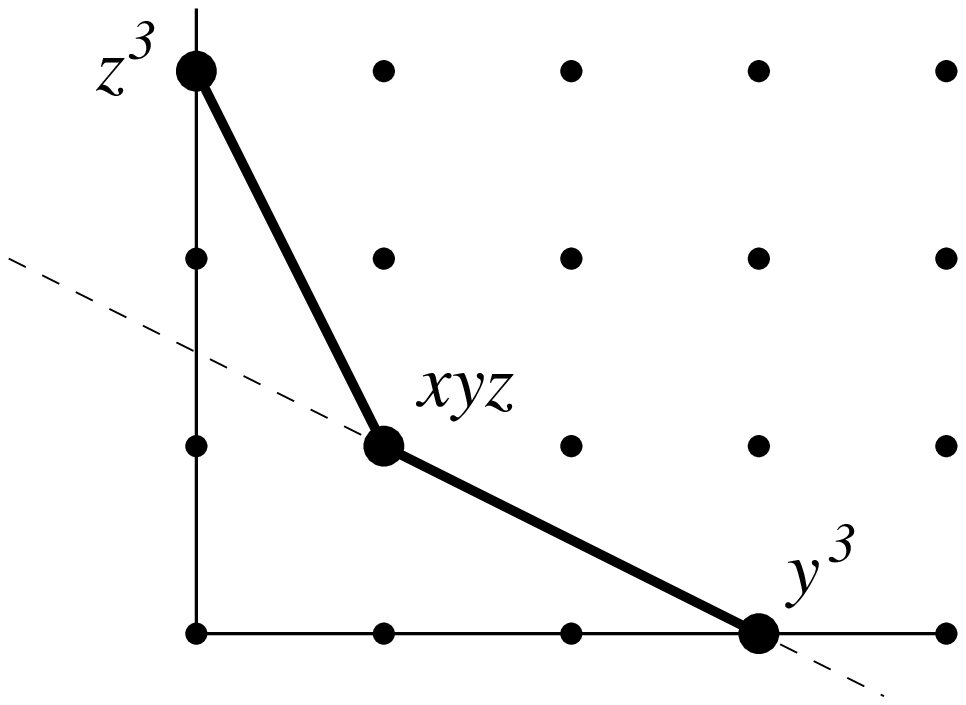}
\end{center}
Setting to zero the coefficient of $z^3$ produces the limit
$y(y^2+xz)$.\vskip 6pt

{\em Type V:\/} Let $p=\im\alpha(0)$ be a singular point of the support
of $\cC$, and let $m$ be the multiplicity of $\cC$ at
$p$. Again choose a line in the tangent cone to $\cC$ at $p$,
and choose coordinates $(x:y:z)$ so that
$p=(1:0:0)$ and $z=0$ is the selected line. 

We may describe $\cC$ near $p$ as the union of $m$ `formal
branches', cf.~\S\ref{formalbranches}; those that are tangent to the line
$z=0$ (but not equal to it) may be written 
$$z=f(y)=\sum_{i\ge 0} \gamma_{\lambda_i} y^{\lambda_i}$$
with $\lambda_i\in \Qbb$, $1<\lambda_0<\lambda_1<\dots$, and
$\gamma_{\lambda_0}\ne 0$.

The choices made above determine a finite set of rational numbers, 
which we call the `characteristics' for $\cC$ (w.r.t.~$p$ and the line $z=0$):
these are the numbers $C$ for which there exist two branches $\cB$, 
$\cB'$ tangent to $z=0$ that agree modulo $y^C$, differ at $y^C$, 
and have $\lambda_0<C$. 
(Formal branches are called `pro-branches' in \cite{MR2107253},
Chapter~4; the numbers~$C$ are `exponents of contact'.)

Let $S$ be the number of branches that agree with $\cB$ (and $\cB'$)
modulo $y^C$. The initial exponents $\lambda_0$ and the
coefficients $\gamma_{\lambda_0}$, $\gamma_{\frac {\lambda_0+C}2}$
for these $S$ branches agree. Let
$\gamma_C^{(1)},\dots,\gamma_C^{(S)}$ be the coefficients of $y^C$
in these branches (so that at least two of these numbers are
distinct, by the choice of~$C$). Then the limit
 is defined by 
$$x^{d-2S}\prod_{i=1}^S\left(zx-\frac {\lambda_0(\lambda_0-1)}2
\gamma_{\lambda_0}y^2 -\frac{\lambda_0+C}2
\gamma_{\frac{\lambda_0+C}2}yx-\gamma_C^{(i)}x^2\right)\quad.$$
This is a union of quadritangent conics with (possibly) a multiple of the
distinguished tangent, which must be supported on the kernel line.

\subsection{The main theorem, and the structure of its proof}\label{proof}
Simple dimension counts show that, for each type as listed in \S\ref{germlist},
the union of the orbits of the marker centers is a set of dimension $7$ in
$\TPbb^8\subset \Pbb^8\times\Pbb^N$; hence it is a dense set in a 
component of the~PNC. In fact, marker centers of type~I, III, IV, and~V have
7-dimensional orbit, so the corresponding components of the PNC are
the orbit closures of these points.

Type~II marker centers are points 
$(\alpha, \cX)\in \Pbb^8\times\Pbb^N$,
where $\alpha$ is a rank-1 matrix whose image is a general point 
of a nonlinear component of $\cC$. The support of~$\cX$ contains a 
conic tangent to the kernel line; this gives
a 1-parameter family of 6-dimensional orbits in $\Pbb^8\times\Pbb^N$,
accounting for a component of the PNC.

We can now formulate a more precise version of Theorem~\ref{main}:
\begin{theorem}[Main theorem]\label{mainmain}
Let $\cC\subset \Pbb^2_\Cbb$ be an arbitrary plane curve. 
The marker germs listed in~\S\ref{germlist} determine components of the 
PNC for $\cC$, as explained above. Conversely, all components of the
PNC are determined by the marker germs of type I--V listed in~\S\ref{germlist}.
\end{theorem}

By the considerations in \S\ref{ident}, this statement implies 
Theorem~\ref{main}.

The first part of Theorem~\ref{mainmain} has been established above.
In order to prove the second part, we will define a simple notion
of `equivalence' of germs 
(Definition~\ref{equivgermsnew}),
such that, in particular, equivalent germs $\alpha(t)$ lead to the same
component of the PNC. We will show that any given germ $\alpha(t)$
centered at a point of $\cS$ either is equivalent (after a parameter change, 
if necessary)
to one of the marker germs, or its lift in~$\TPbb^8$ meets the PNC at a point
of $R$  (cf.~Remark~\ref{eluding}) or of the boundary of the orbit of a marker center.
In the latter cases, the center of the lift varies in a locus of 
dimension~$<7$, hence such germs do not contribute components to the PNC. 
The following lemma allows us to identify easily limits in the intersection
of $R$ and the PNC.

\begin{lemma}\label{rank2lemma}
Assume that $\alpha(0)$ has rank~$1$.
If $\lim_{t\to 0}\cC\circ\alpha(t)$ is a star with center on 
$\ker\alpha(0)$, then it is a rank-2 limit.
\end{lemma}

\begin{proof}
Assume $\cX=\lim_{t\to 0}\cC\circ\alpha(t)$ is a
star with center on $\ker\alpha(0)$. We may choose coordinates so that
$x=0$ is the kernel line, and the generator for the ideal of
$\cX$ is a polynomial in $x,y$ only. If
$$\alpha(t)=\begin{pmatrix}
a_{11}(t) & a_{12}(t) & a_{13}(t) \\
a_{21}(t) & a_{22}(t) & a_{23}(t) \\
a_{31}(t) & a_{32}(t) & a_{33}(t)
\end{pmatrix}\quad,$$
then $\cX=\lim_{t\to 0}\cC\circ\beta(t)$ for
$$\beta(t)=\begin{pmatrix}
a_{11}(t) & a_{12}(t) & 0 \\
a_{21}(t) & a_{22}(t) & 0 \\
a_{31}(t) & a_{32}(t) & 0
\end{pmatrix}\quad.$$
Since $\alpha(0)$ has rank~1 and kernel line $x=0$,
$$\alpha(0)=\begin{pmatrix}
a_{11}(0) & 0 & 0 \\
a_{21}(0) & 0 & 0 \\
a_{31}(0) & 0 & 0
\end{pmatrix}=\beta(0)\quad.$$
Now $\beta(t)$ is contained in the rank-2 locus, verifying the assertion.
\end{proof}

A limit $\lim_{t\to 0}\cC\circ\alpha(t)$ as in this lemma will be called a 
`kernel star'.

Sections~\ref{setth} and~\ref{typeVcomps} contain the successive reductions 
bringing a given 
germ $\alpha(t)$ centered at a point of $\cS$ into one of the forms given in 
\S\ref{germlist}, or establishing that it does not contribute a component of 
the PNC. This analysis will conclude the proof of Theorem~\ref{mainmain}.

\subsection{Two examples}\label{twoexamples}
The two examples that follow illustrate the main result, and 
show that components of all types may already occur on curves of degree~4.
Simple translations are used to bring the marker germs provided by 
\S\ref{germlist} into the form given here.

\begin{example}
Consider the reducible quartic~$\cC_1$ given by the equation
$$(y+z)(xy^2+xyz+xz^2+y^2z+yz^2)=0\,.$$
It consists of an irreducible cubic with a node at $(1:0:0)$ and a line through
the node and the inflection point $(0:1:-1)$. The other inflection points
are $(0:1:0)$ and $(0:0:1)$. 
According to Theorem~\ref{mainmain} and the list in \S\ref{germlist},
the PNC for~$\cC_1$ has one component of
type I, one component of type~II, one component of type III, corresponding
to the triple point $(1:0:0)$, and four components of type IV: one
for each of the inflection points $(0:1:0)$ and $(0:0:1)$, one for the
node $(0:1:-1)$ and the tangent line $x=y+z$ to the cubic at that point,
and one for the triple point $(1:0:0)$ and the two lines
in the tangent cone $y^2+yz+z^2=0$ to the cubic at that point.
Here is a schematic drawing of the curve, with features marked by the
corresponding types (four points are marked as $\text{IV}_i$, since
four different points are responsible for the presence of type IV components):
\begin{center}
\includegraphics[scale=.45]{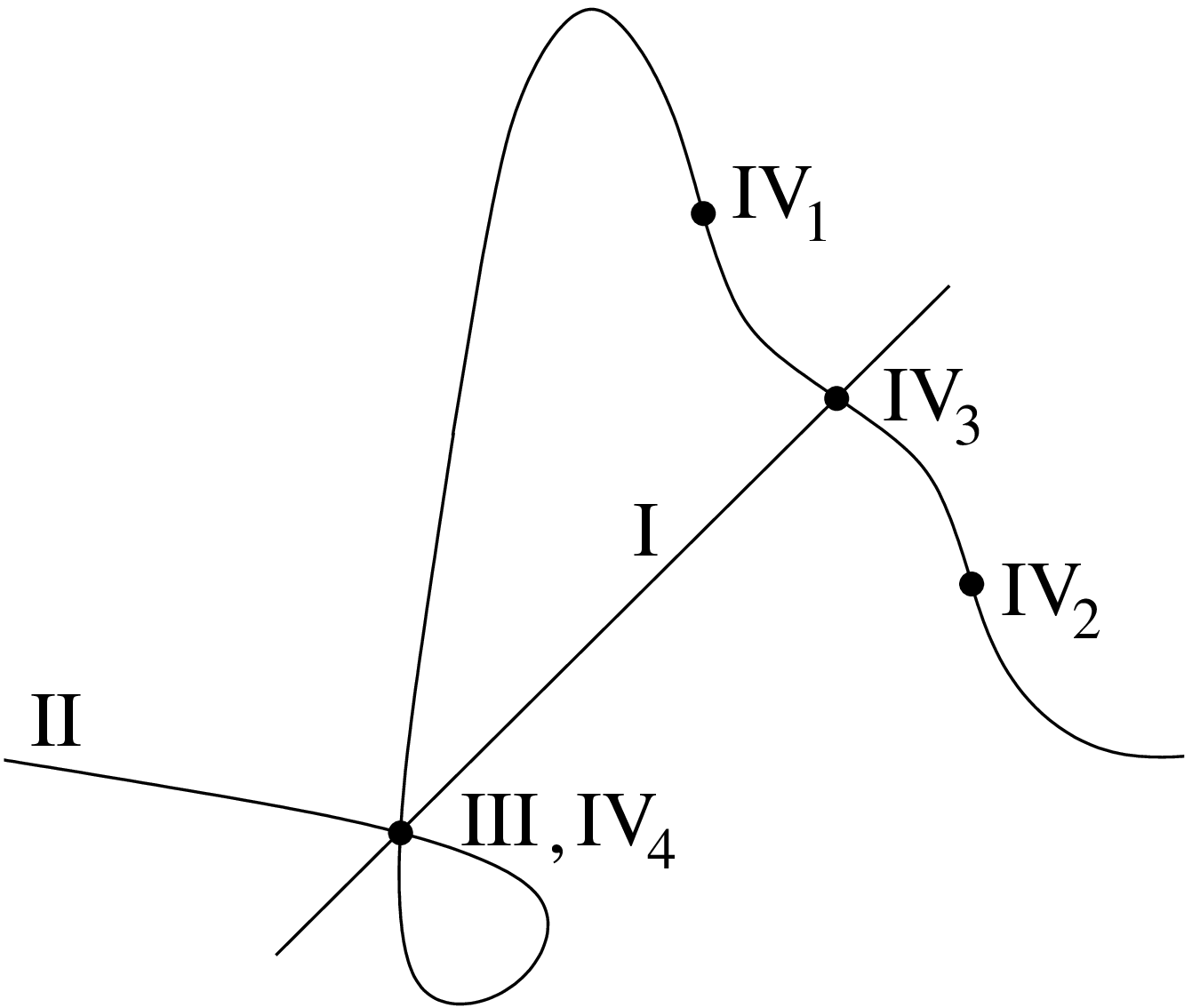}
\end{center}
A list of representative marker germs is as follows:
$$
{\rm I}:
\begin{pmatrix}
1 & 0 & 0\\
0 & 1 & 0\\
0 & -1 & t
\end{pmatrix};\quad
{\rm II}:
\begin{pmatrix}
2 & 0 & 0\\
-3 & t & 0\\
6 & 0 & t^2
\end{pmatrix};\quad
{\rm III}:
\begin{pmatrix}
1 & 0 & 0\\
0 & t & 0\\
0 & 0 & t
\end{pmatrix};
$$
and, for type~IV:
$$
\begin{pmatrix}
t & 0 & 0\\
0 & 1 & 0\\
-t & 0 & t^3
\end{pmatrix},\quad
\begin{pmatrix}
t & 0 & 0\\
-t & t^3 & 0\\
0 & 0 & 1
\end{pmatrix},\quad
\begin{pmatrix}
t & 0 & 0\\
0 & 1 & 0\\
t & -1 & t^3
\end{pmatrix},\quad
\begin{pmatrix}
1 & 0 & 0\\
0 & \rho t & 0\\
0 & t & t^2
\end{pmatrix},\quad
\begin{pmatrix}
1 & 0 & 0\\
0 & \rho^2 t & 0\\
0 & t & t^2
\end{pmatrix}
$$
(where $\rho$ is a primitive third root of unity).
The latter two marker germs have the same center and lead to 
projectively equivalent limits, hence they contribute the same component 
of the PNC. 
The corresponding limits of $\cC_1$ are given by
$$
xy^2z,\quad
x^2(8y^2-9xz),\quad
x(y+z)(y^2+yz+z^2),\quad
y(y^2z+x^3),\quad
z(yz^2+x^3),
$$$$
x(y^2z-x^3),\quad
y^2(y^2-(\rho+2)xz),\quad \text{and} \quad
y^2(y^2-(\rho^2+2)xz),
$$
respectively: a triangle with one line doubled,
a conic with a double tangent line,
a fan with star centered at $(1:0:0)$,
a cuspidal cubic with its cuspidal tangent
(two limits),
a cuspidal cubic with the line through the cusp and the inflection point,
and finally  
a conic with a double transversal line
(two limits). Schematically, the limits may be represented as follows: 
\begin{center}
\includegraphics[scale=.45]{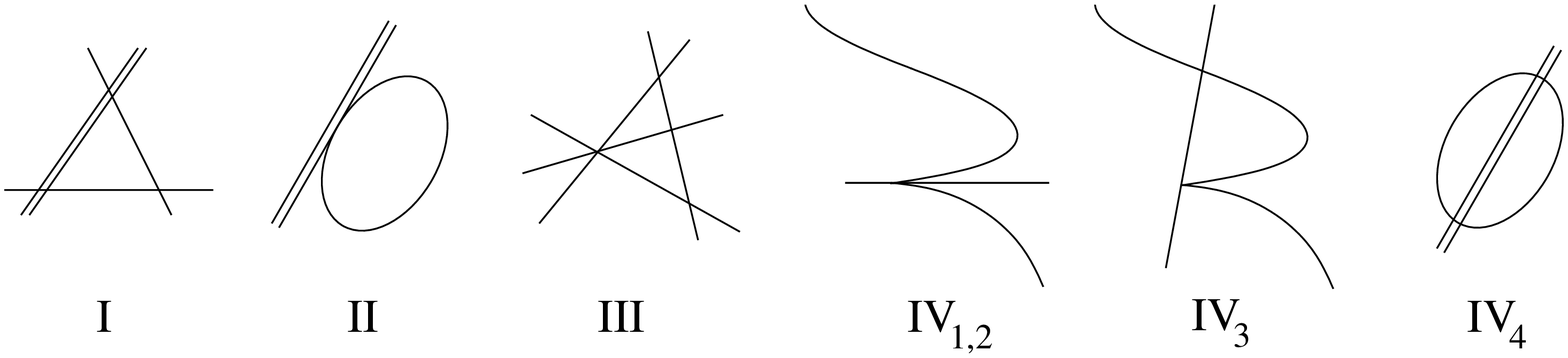}
\end{center}
According to Theorem~\ref{main}, all limits of $\cC_1$ (other than stars of
lines) are projectively equivalent to one of these curves, or to limits of them
(cf.~\S\ref{boundary}).
\qed\end{example}

\begin{example}\label{extwo}
Consider the irreducible quartic~$\cC_2$ given by the equation
$$(y^2-xz)^2=y^3z.$$
It has a ramphoid cusp at $(1:0:0)$, an ordinary cusp at $(0:0:1)$, and
an ordinary inflection point at $(3^3 5{:}{-}2^6 3^2{:}{-}2^{12})$; there
are no other singular or inflection points.
The PNC for~$\cC_2$ has one component of type~II, two components of type~IV,
corresponding to the inflection point and the ordinary cusp, and one
component of~type~V, corresponding to the ramphoid cusp. (Note that there is
no component of type~IV corresponding to the ramphoid cusp.)
Representative marker germs for the latter two components are
$$
{\rm IV}:
\begin{pmatrix}
0 & t^3 & 0\\
t^2 & 0 & 0\\
0 & 0 & 1
\end{pmatrix}\quad{\rm and}\quad
{\rm V}:
\begin{pmatrix}
1 & 0 & 0\\
t^4 & t^5 & 0\\
t^8 & 2t^9 & t^{10}
\end{pmatrix}
$$
and the corresponding limits of $\cC_2$ are given by
$$z(y^2z-x^3)\quad{\rm and}\quad(y^2-xz+x^2)(y^2-xz-x^2),$$
respectively: a cuspidal cubic with its inflectional tangent and 
a pair of quadritangent conics. The connected component of the
stabilizer of the latter limit is the additive group. The germ with
entries $1$, $t$, and $t^2$ on the diagonal and zeroes elsewhere
leads to the limit $(y^2-xz)^2$, a double conic; its orbit is too
small to produce an additional component of type IV.
\qed\end{example}

\section{Proof of the main theorem: key reductions and components of
type~I--IV}\label{setth}

\subsection{Outline}\label{preamble}
In this section we show that, for a given curve $\cC$, any germ $\alpha(t)$
contributing to the PNC is `equivalent'
(up to a coordinate and parameter change, if necessary)
to a marker germ as listed in \S\ref{germlist}. 
As follows from~\S\ref{rough} and Lemma~\ref{PNCtolimits},
we may assume that $\det\alpha(t)\not\equiv 0$ and that the image 
of~$\alpha(0)$ is contained in~$\cC$.

Observe that if the center $\alpha(0)$ has rank~2 and is
a point of $\cS$, then $\alpha(t)$ is already of the form given in
\S\ref{germlist}, Type~I; it is easy to verify that the limit is then
as stated there. This determines completely the components of type~I.
Thus, we will assume in most of what follows that $\alpha(0)$ has 
rank~1, and its image is a point of~$\cC$.

\subsubsection{Equivalence of germs}
\begin{defin}\label{equivgermsnew}
Two germs $\alpha(t)$, $\beta(t)$ are {\em equivalent\/} if 
$\beta(t\nu(t))\equiv \alpha(t)\circ m(t)$, with $\nu(t)$ a unit
in $\Cbb[[t]]$, and $m(t)$ a germ such that $m(0)=I$ (the identity).
\end{defin}

For example: if $n(t)$ is a $\Cbb[[t]]$-valued point of $\PGL(3)$,
then $\alpha(t)\circ n(t)$ is equivalent to $\alpha(t)\circ n(0)$. We
will frequently encounter this situation.

\begin{lemma}\label{stequivnew}
Let $\cC$ be any plane curve, with defining homogeneous ideal
$(F(x,y,z))$. If $\alpha(t)$, $\beta(t)$ are equivalent germs, then
the initial terms in $F\circ\alpha(t)$, $F\circ\beta(t)$ coincide up
to a nonzero multiplicative constant;
in particular, the limits $\lim_{t\to 0}\cC\circ \alpha(t)$, 
$\lim_{t\to 0}\cC\circ \beta(t)$ are equal.\hfill\qed
\end{lemma}

If $\alpha$ and $\beta$ are equivalent germs, note that $\alpha(0)=
\beta(0)$; by Lemma~\ref{stequivnew} it follows that, for every curve 
$\cC$, $\alpha$ and $\beta$ lift to germs in $\TPbb^8$ centered at
the same point.

\subsubsection{Summary of the argument}
The general plan for the rest of this section  is as follows: we will show
that every contributing $\alpha(t)$ centered at a rank-1 matrix is equivalent
(in suitable coordinates, and possibly up to a parameter change)
to one of the form
$$\alpha(t)=\begin{pmatrix}
1 & 0 & 0\\
0 & t^b & 0\\
0 & 0 & t^c
\end{pmatrix}\quad\text{or}\quad
\begin{pmatrix}
1 & 0 & 0     \\
t^a & t^b & 0 \\
\underline{f(t^a)} & \underline{f'(t^a) t^b} & t^c\end{pmatrix}
\quad,$$
where $b\le c$ resp.~$a<b\le c$ are positive integers,
$z=f(y)$ is a formal branch for $\cC$ at $(1:0:0)$, and $\underline{\cdots}$ 
denotes the truncation modulo $t^c$ (cf.~\S\ref{germlist} and \S\ref{details}).

The main theorem will follow from further analyses of these forms, identifying
which do {\em not\/} contribute components to the PNC, and leading
to the restrictions explained in \S\ref{germlist} and \S\ref{details}.
Specifically, the germs on the left lead to components of type~II, III, and IV
(\S\ref{1PS}); those on the right lead to components of type~V. The latter
germs require a subtle study, performed
in \S\ref{typeVcomps}, leading to the definition of 
`characteristics' and to the description given in \S\ref{details}
(cf.~Proposition~\ref{typeV}).

\subsection{Linear algebra}

\subsubsection{}
This subsection is devoted to the proof of the following result.

\begin{prop}\label{keyreduction}
Every germ as specified in \S\ref{preamble} is equivalent to 
one which, up to a parameter change, has matrix representation
$$\begin{pmatrix}
1 & 0 & 0\\
q(t) & t^b & 0\\
r(t) & s(t)t^b & t^c
\end{pmatrix}$$
in suitable coordinates, with $1\le b\le c$ and $q,r,s$ polynomials
such that $\deg(q)<b$, $\deg(r)<c$, $\deg(s)<c-b$, and
$q(0)=r(0)=s(0)=0$.
\end{prop}
A refined version of this statement is given in Lemma~\ref{faber}.

We will deal with $3\times 3$ matrices with
entries in $\Cbb[[t]]$, that is, $\Cbb[[t]]$-valued points of $\Hom(V,W)$,
for $V$, $W$ 3-dimensional complex vector spaces with chosen bases.
Every such matrix $\alpha(t)$ determines a germ in $\Pbb^8$.
A generator $F$ of the ideal of $\cC$ will be viewed as an element of 
$\Sym^d W^*$, for $d=\deg\cC$; the composition $F\circ\alpha(t)$, a
$\Cbb[[t]]$-valued point of $\Sym^d V^*$, generates the ideal of 
$\cC\circ\alpha(t)$.

We will call matrices of the form
$$\lambda(t)=\begin{pmatrix}
t^a & 0 & 0\\
0 & t^b & 0\\
0 & 0 & t^c
\end{pmatrix}$$
`1-PS', as 
they correspond to 1-parameter subgroups of $\PGL(3)$.

We will say that two matrices $\alpha(t)$, $\beta(t)$ are equivalent
if the corresponding germs are equivalent in the sense of 
Definition~\ref{equivgermsnew}.
The following lemma will allow us to simplify matrix expressions of germs 
up to equivalence. Define the degree of the zero polynomial to be $-\infty$. 

\begin{lemma}\label{MPI} Let
$$ h_1(t) =\begin{pmatrix}
u_1 & b_1 & c_1 \\
a_2 & u_2 & c_2 \\
a_3 & b_3 & u_3
\end{pmatrix}
$$
be a matrix with entries in $\Cbb[[t]]$, such that $h_1(0)=I$,
and let $a\le b\le c$ be integers.
Then $h_1(t)$ can be written as a product $h_1(t)=h(t)\cdot j(t)$, with
$$
h(t)=\begin{pmatrix}
1 & 0 & 0 \\
q & 1 & 0 \\
r & s & 1 
\end{pmatrix}
\quad,\quad j(t)=\begin{pmatrix}
v_1 & e_1 & f_1 \\
d_2 & v_2 & f_2 \\
d_3 & e_3 & v_3 
\end{pmatrix}
$$
where $q$, $r$, $s$ are {\em polynomials,\/} satisfying
\begin{enumerate}
\item $h(0)=j(0)=I$;
\item $\deg(q)<b-a$, $\deg(r)<c-a$, $\deg(s)<c-b$;
\item\label{refMPI} $d_2\equiv0\pmod{t^{b-a}}$, $d_3\equiv0\pmod{t^{c-a}}$, 
$e_3\equiv0\pmod{t^{c-b}}$.
\end{enumerate}
\end{lemma}

\begin{proof} Necessarily $v_1=u_1, e_1=b_1$ and $f_1=c_1$. Use division 
with remainder to write 
$ v_1^{-1}a_2=D_2t^{b-a}+q $
with $\deg(q)<b-a$, and let $d_2=v_1D_2t^{b-a}$ (so that $qv_1+d_2=a_2$).
This defines $q$ and $d_2$, and uniquely determines $v_2$ and $f_2$.
(Note that $q(0)=d_2(0)=f_2(0)=0$ and that $v_2(0)=1$.)

Similarly, we let $r$ be the remainder of
$(v_1v_2-e_1d_2)^{-1}(v_2a_3-d_2b_3)$
after division by $t^{c-a}$; and $s$ be the remainder of
$(v_1v_2-e_1d_2)^{-1}(v_1b_3-e_1a_3)$
after division by $t^{c-b}$.
Then $\deg(r)<c-a$, $\deg(s)<c-b$ and $r(0)=s(0)=0$; moreover, we have
$$ v_1r+d_2s\equiv a_3\pmod{t^{c-a}},\qquad e_1r+v_2s\equiv b_3
\pmod{t^{c-b}},$$
so we take $d_3=a_3-v_1r-d_2s$ and $e_3=b_3-e_1r-v_2s$. This defines $r$,
$s$, $d_3$ and $e_3$, and uniquely determines $v_3$.
\end{proof}

\begin{corol}\label{MPIcorol}
Let $h_1(t)$ be a matrix with entries in $\Cbb[[t]]$, such that $h_1(0)=I$,
and let $a\le b\le c$ be integers. Then there exists a constant invertible
matrix $L$ such that the product
$$
h_1(t)\cdot
\begin{pmatrix}
t^a & 0 & 0 \\
0 & t^b & 0 \\
0 & 0 & t^c
\end{pmatrix}
$$
is equivalent to
$$
\begin{pmatrix}
1 & 0 & 0 \\
q & 1 & 0 \\
r & s & 1 
\end{pmatrix}\cdot
\begin{pmatrix}
t^a & 0 & 0 \\
0 & t^b & 0 \\
0 & 0 & t^c
\end{pmatrix}\cdot L
$$
where $q$, $r$, $s$ are polynomials such that $\deg(q)<b-a$, $\deg(r)<c-a$, 
$\deg(s)<c-b$, and $q(0)=r(0)=s(0)=0$.
\end{corol}

\begin{proof}
With notation as in Lemma~\ref{MPI} we have 
$$j(t)\cdot 
\begin{pmatrix}
t^a & 0 & 0 \\
0 & t^b & 0 \\
0 & 0 & t^c
\end{pmatrix}
=
\begin{pmatrix}
v_1 t^a & e_1 t^b & f_1 t^c \\
d_2 t^a & v_2 t^b & f_2 t^c \\
d_3 t^a & e_3 t^b & v_3 t^c
\end{pmatrix}
=
\begin{pmatrix}
t^a & 0 & 0 \\
0 & t^b & 0 \\
0 & 0 & t^c
\end{pmatrix}
\cdot \ell(t)\quad,$$
with
$$\ell(t)=\begin{pmatrix}
v_1 & e_1 t^{b-a} & f_1 t^{c-a} \\
d_2 t^{a-b} & v_2  & f_2 t^{c-b} \\
d_3 t^{a-c} & e_3 t^{b-c} & v_3
\end{pmatrix}\quad.$$
By (\ref{refMPI}) in Lemma~\ref{MPI}, $\ell(t)$ has entries in $\Cbb[[t]]$ and
is invertible; in fact, $L=\ell(0)$ is lower triangular, 
with 1's on the diagonal. Therefore Lemma~\ref{MPI} gives
$$
h_1(t)\cdot
\begin{pmatrix}
t^a & 0 & 0 \\
0 & t^b & 0 \\
0 & 0 & t^c
\end{pmatrix}=
\begin{pmatrix}
1 & 0 & 0 \\
q & 1 & 0 \\
r & s & 1 
\end{pmatrix}\cdot
\begin{pmatrix}
t^a & 0 & 0 \\
0 & t^b & 0 \\
0 & 0 & t^c
\end{pmatrix}\cdot \ell(t)\quad,
$$
from which the statement follows.
\end{proof}

The gist of this result is that, up to equivalence, matrices `to the left of a 
1-PS' and centered at the identity may be assumed to be lower triangular,  and 
to have polynomial entries, with controlled degrees.

\subsubsection{}
We denote by $v$ the order of vanishing at~$0$
of a polynomial or power series; we define $v(0)$ to be $+\infty$.
The following statement is a refined version of Proposition~\ref{keyreduction}.

\begin{lemma}\label{faber}
Let $\alpha(t)$ be a $\/3\times 3$ matrix with entries in $\Cbb[[t]]$,
such that $\alpha(0)\ne 0$ and $\det \alpha(t)\not\equiv 0$.
Then there exist constant invertible matrices $H$, $M$ such that
$\alpha(t)$ is equivalent to
$$\beta(t)=H\cdot \begin{pmatrix}
1 & 0 & 0 \\
q & 1 & 0 \\
r & s & 1
\end{pmatrix} \cdot \begin{pmatrix}
1 & 0 & 0 \\
0 & t^b & 0 \\
0 & 0 & t^c
\end{pmatrix} \cdot M\quad,$$
with
\begin{itemize}
\item $b\le c$ nonnegative integers, $q,r,s$ polynomials;
\item $\deg(q)<b$, $\deg(r)<c$, $\deg(s)<c-b$;
\item $q(0)=r(0)=s(0)=0$.
\end{itemize}
If, further, $b=c$ and $q$, $r$ are not both zero, then we may assume that
$v(q)<v(r)$. 

Finally, if $q(t)\not\equiv 0$ then we may choose $q(t)=t^a$, with $a=v(q)<b$
(and thus $a<v(r)$ if $b=c$).
\end{lemma}

\begin{proof}
By standard diagonalization of matrices over Euclidean domains,
every $\alpha(t)$ as in the statement can be written as a product
$$h_0(t)\cdot
\begin{pmatrix}
1 & 0 & 0 \\
0 & t^b & 0 \\
0 & 0 & t^c
\end{pmatrix}
\cdot k(t)\quad,$$
where $b\le c$ are nonnegative integers, and $h_0(t)$, $k(t)$ are invertible 
(over $\Cbb[[t]]$). Letting $H=h_0(0)$, $h_1(t)=H^{-1}\cdot h_0(t)$, and
$K=k(0)$, this shows that $\alpha(t)$ is equivalent to
$$H\cdot h_1(t) \cdot 
\begin{pmatrix}
1 & 0 & 0 \\
0 & t^b & 0 \\
0 & 0 & t^c
\end{pmatrix}
\cdot K$$
with $h_1(0)=I$, and $K$ constant and invertible.
By Corollary~\ref{MPIcorol}, this matrix is equivalent to
$$\beta(t)=H\cdot \begin{pmatrix}
1 & 0 & 0 \\
q & 1 & 0 \\
r & s & 1 \end{pmatrix}\cdot \begin{pmatrix}
1 & 0 & 0 \\
0 & t^b & 0 \\
0 & 0 & t^c
\end{pmatrix}
\cdot L\cdot K$$
with $L$ invertible, and $q$, $r$, $s$ polynomials satisfying the needed
conditions. Letting $M=L\cdot K$ gives the statement in the case $b<c$.

If $b=c$, then the condition that $\deg s<c-b=0$ forces
$s\equiv 0$. When $q$ and $r$ are not both~$0$, the inequality
$v(q)<v(r)$ may be obtained by conjugating with a constant matrix.

If $q(t)\not\equiv 0$ and $v(q)=a$, then we can extract its $a$-th root as a 
power series. It follows that there exists a unit $\nu(t)\in\Cbb[[t]]$ such that
$q(t\nu(t))=t^a$.
Therefore,
$$\beta(t\nu(t))=H\cdot \begin{pmatrix}
1 & 0 & 0 \\
t^a & 1 & 0 \\
r(t\nu(t)) & s(t\nu(t)) & 1 \end{pmatrix}\cdot \begin{pmatrix}
1 & 0 & 0 \\
0 & t^b & 0 \\
0 & 0 & t^c
\end{pmatrix}
\cdot \begin{pmatrix}
1 & 0 & 0 \\
0 & \nu(t)^b & 0 \\
0 & 0 & \nu(t)^c
\end{pmatrix}
\cdot M\quad.$$
Another application of Corollary~\ref{MPIcorol} allows us to truncate the
power series $r(t\nu(t))$ and $s(t\nu(t))$ to obtain polynomials $\underline r$,
$\underline s$ satisfying the same conditions as $r$, $s$, at the price
of multiplying to the right of the 1-PS by a constant invertible matrix 
$\underline K$: that is, $\beta(t\nu(t))$ (and hence $\alpha(t)$) is equivalent 
to
$$H\cdot \begin{pmatrix}
1 & 0 & 0 \\
t^a & 1 & 0 \\
\underline r & \underline s & 1 \end{pmatrix}\cdot \begin{pmatrix}
1 & 0 & 0 \\
0 & t^b & 0 \\
0 & 0 & t^c
\end{pmatrix}
\cdot \left[ \underline K \cdot \begin{pmatrix}
1 & 0 & 0 \\
0 & \nu(0)^b & 0 \\
0 & 0 & \nu(0)^c
\end{pmatrix} 
\cdot M \right]\quad.$$
Renaming $r=\underline r$, $s=\underline s$, and absorbing the factors on
the right into $M$ completes the proof of Lemma~\ref{faber}.
\end{proof}

The matrices $H$, $M$ appearing in Lemma~\ref{faber} may be
omitted by changing the bases of $W$ and $V$ accordingly. 
Further, we may assume that $b>0$, since we are
already reduced to the case in which $\alpha(0)$ is a rank-1
matrix. This concludes the proof of Proposition~\ref{keyreduction}.
In what follows, we will assume that $\alpha$ is a germ in the standard 
form given above.

\subsection{Components of type II, III, and IV}\label{1PS}
It will now be convenient to switch to affine
coordinates centered at the point $(1:0:0)$. We write
$$F(1:y:z)=F_m(y,z)+F_{m+1}(y,z)+\cdots +F_d(y,z)\quad,$$
with $d=\deg \cC$, $F_i$ homogeneous of degree $i$, and $F_m\ne
0$. Thus, $F_m(y,z)$ generates the ideal of the {\em tangent cone\/}
of $\cC$ at $p$.

We first consider the case in which $q=r=s=0$,
that is, in which $\alpha(t)$ is itself a 1-PS:
$$\alpha(t)=\begin{pmatrix}
1 & 0 & 0\\
0 & t^b & 0\\
0 & 0 & t^c
\end{pmatrix}$$
with $1\le b \le c$. Also, we may assume that $b$ and $c$ are coprime:
this only amounts to a reparametrization of the germ by $t \mapsto
t^{1/gcd(b,c)}$; the new germ is not equivalent to the old
one in terms of Definition~\ref{equivgermsnew}, but clearly achieves the
same limit.

Germs with $b=c$ $(=1)$ lead to components of type~III, cf.~\S\ref{germlist}
(also cf.~\cite{MR2001h:14068}, \S2, Fact~4(i)):

\begin{prop}\label{tgcone}
If $q=r=s=0$ and $b=c$, then $\lim_{t\to 0} \cC\circ\alpha(t)$
is a fan consisting of a star projectively equivalent to the tangent cone
to $\cC$ at $p$, and of a residual $(d-m)$-fold line supported
on $\ker\alpha$.
\end{prop}

\begin{proof}
The composition $F\circ\alpha(t)$ is
$$F(x:t^by:t^bz)=t^{bm}x^{d-m}F_m(y,z)+t^{b(m+1)}x^{d-(m+1)}
F_{m+1}(y,z)+\cdots+ t^{dm}F_d(y,z)\quad.$$
By definition of limit, $\lim_{t\to 0}\cC\circ\alpha(t)$ has ideal
$(x^{d-m}F_m(y,z))$, proving the assertion.
\end{proof}

The case $b<c$ corresponds to the germs of type~II and type~IV 
in \S\ref{germlist}. We have to prove that contributing germs of this
type are precisely those satisfying the further restrictions specified
there: specifically, $-b/c$ must be a slope of one of the Newton polygons 
for $\cC$ at the point. 
We first show that $z=0$ must be a component of the tangent cone:

\begin{lemma}
If $q=r=s=0$ and $b<c$, and $z=0$ is not contained in the tangent cone
to $\cC$ at $p$, then $\lim_{t\to 0} \cC\circ\alpha(t)$
is a rank-2 limit.
\end{lemma}

\begin{proof}
The condition regarding $z=0$ translates into $F_m(1,0)\ne 0$. 
Applying $\alpha(t)$ to~$F$, we find:
$$F(x:t^by:t^cz)=t^{bm}x^{d-m} F_m(y,t^{c-b}z)+t^{b(m+1)} x^{d-(m+1)}
F_{m+1}(y,t^{c-b}z)+\cdots$$
Since $F_m(1,0)\ne 0$, the dominant term on the right-hand-side is
$x^{d-m}y^m$. This proves the assertion, by Lemma~\ref{rank2lemma}.
\end{proof}

Components of the PNC that arise due to 1-PS with $b<c$ may be 
described in terms of the {\em Newton polygon\/} for $\cC$ at $(0,0)$
relative to the line $z=0$, which we may now assume 
to be part of the tangent cone to $\cC$ at $p$. The Newton
polygon for $\cC$ in the chosen coordinates is the boundary of
the convex hull of the union of the positive quadrants with origin at
the points $(j,k)$ for which the coefficient of $x^iy^jz^k$ in the
equation for $\cC$ is nonzero (see \cite{MR88a:14001},
p.~380). 
The part of the Newton polygon consisting of line segments
with slope strictly between $-1$ and $0$ does not depend on the choice
of coordinates fixing the flag $z=0$, $p=(0,0)$.

\begin{prop}\label{Newtonsides}
Assume $q=r=s=0$ and $b<c$.
\begin{itemize}
\item If $-b/c$ is not a slope of the Newton polygon for $\cC$,
  then the limit $\lim_{t\to 0} \cC\circ\alpha(t)$ is supported
  on (at most) three lines; these curves do not contribute 
  components to the PNC.
\item If $-b/c$ is a slope of a side of the Newton polygon for
  $\cC$, then the ideal of the limit $\lim_{t\to 0}\cC
  \circ\alpha(t)$ is generated by the polynomial obtained by 
  setting to~$0$ the coefficients of
  the monomials in $F$ {\em not\/} on that side.
  Such polynomials are of the form
$$G=x^{\overline e}y^fz^e \prod_{j=1}^S(y^c+\rho_j x^{c-b}z^b)
\quad.$$
\end{itemize}
\end{prop}

\begin{proof}
For the first assertion, simply note that under the stated hypotheses
only one monomial in $F$ is dominant in $F\circ\alpha(t)$; hence, the
limit is supported on the union of the coordinate axes. A simple
dimension count shows that such limits may span at most a 6-dimensional
locus in $\Pbb^8\times\Pbb^N$, and it follows that such germs do
not contribute a component to the PNC.

For the second assertion, note that the dominant terms in
$F\circ\alpha(t)$ are precisely those on the side of the Newton polygon
with slope equal to $-b/c$. It is immediate that the resulting
polynomial can be factored as stated.
\end{proof}

If the point $p=(1:0:0)$ is a singular or an inflection point of
the support of~$\cC$, and $b/c\ne 1/2$, we find the type~IV
germs of \S\ref{germlist}; also cf.~\cite{MR2001h:14068}, \S2, 
Fact~4(ii). The number $S$ of `cuspidal' factors in
$G$ is the number of segments cut out by the integer lattice on the
selected side of the Newton polygon.
If $b/c=1/2$, then a dimension count shows that the corresponding limit 
will contribute a component to the PNC (of type~IV)
unless it is supported on a conic union (possibly) the kernel line.

If $p$ is a {\em nonsingular, non-inflectional\/} point of the support
of $\cC$, 
then the Newton polygon consists of a single side
with slope $-1/2$; these are the type~II germs of \S\ref{germlist}.
Also cf.~\cite{MR2001h:14068}, Fact~2(ii).

\section{Components of type~V}\label{typeVcomps}
Having dealt with the 1-PS case in the previous section, we may now
assume that 
\begin{equation*}
\tag{$\dagger$}
\alpha(t)=\begin{pmatrix}
1 & 0 & 0\\
q(t) & t^b & 0\\
r(t) & s(t)t^b & t^c
\end{pmatrix}
\end{equation*}
with the conditions listed in Lemma~\ref{faber}, and further
{\em such that $q,r$, and $s$ do not all vanish identically.\/}

Our task is to show that contributing germs of this kind must in
fact be of the form specified in \S\ref{germlist} and \S\ref{details}.
We will show that a germ $\alpha(t)$ as above leads to
a rank-2 limit (and hence does not contribute a component to the PNC)
unless $\alpha(t)$ and certain formal branches (cf.~\cite{MR88a:14001} 
and \cite{MR1836037}, Chapter~6 and~7) of the curve are closely related. 
More precisely, we will prove the following result.

\begin{prop}\label{standardform}
Let $\alpha(t)$ be as specified above, and assume that 
$\lim_{t\to 0}\mathcal C\circ\alpha(t)$ is not a rank-2 limit. Then 
$\cC$ has a formal branch $z=f(y)$, tangent to $z=0$, such that 
$\alpha$ is equivalent to a germ 
$$\begin{pmatrix}
1 & 0 & 0     \\
t^a & t^b & 0 \\
\underline{f(t^a)} & \underline{f'(t^a) t^b} & t^c\end{pmatrix}
\quad,$$
with $a<b<c$ positive integers.
Further, it is necessary that $\frac ca\le \lambda_0+2(\frac ba-1)$,
where $\lambda_0>1$ is the (fractional) order of the branch.
\end{prop}

\noindent For a power series $g(t)$ with fractional exponents, we 
write here $\underline{g(t)}$ for its truncation modulo $t^c$. 
(The truncations appearing in the statement are in fact polynomials.)

The proof of the proposition requires the analysis of several
cases. We will first show that under the hypothesis that $\lim_{t\to
  0} \cC\circ\alpha(t)$ is not a rank-2 limit we may assume
that $q(t)\not\equiv 0$, and this will allow us to replace it with a power of
$t$; next, we will deal with the $b=c$ case; and finally we will see
that if $b<c$ and $\alpha(t)$ is not in the stated form, then the
limit of every irreducible branch of $\cC$ is a star with center $(0:0:1)$. 
This will imply that the limit of $\cC$ is a kernel star in this case, proving 
the assertion by Lemma~\ref{rank2lemma}.

This analysis is carried out in \S\ref{formalbranches}--\ref{Eop}.
In \S\ref{charaV} we determine germs of the form given
in Proposition~\ref{standardform} that can lead to components of
type~V, obtaining the description given in \S\ref{germlist}.
In \S\ref{quadritangent} we complete the proof of 
Theorem~\ref{mainmain}, recovering the description given in 
\S\ref{details} of the limits obtained along these germs.

\subsection{Limits of formal branches}\label{formalbranches}
In this subsection we recall the notion of formal branches and define
the `limit' of a formal branch. The limit of a curve $\cC$ will be
expressed in terms of the limits of its formal branches.

Choose affine coordinates $(y,z)=(1:y:z)$ so that $p=(0,0)$,
and let $\Phi(y,z)=F(1:y:z)$ be the generator for the ideal of
$\cC$ in these coordinates. Decompose $\Phi(y,z)$ in
$\Cbb[[y,z]]$:
$$\Phi(y,z)=\Phi_1(y,z)\cdot\cdots\cdot \Phi_r(y,z)$$
with $\Phi_i(y,z)$ irreducible power series. 
These define the {\em irreducible branches\/} of $\cC$
at~$p$. Each $\Phi_i$ has a unique tangent line at $p$; if this
tangent line is {\em not\/} $y=0$, by the Weierstrass preparation
theorem we may write (up to a unit in $\Cbb[[y,z]]$) $\Phi_i$ as a
monic polynomial in $z$ with coefficients in $\Cbb[[y]]$, of degree
equal to the multiplicity $m_i$ of the branch at $p$
(cf.~for example \cite{MR1836037}, \S6.7).
If $\Phi_i$ {\em is\/} tangent to $y=0$, we may likewise write it as a
polynomial in $y$ with coefficients in $\Cbb[[z]]$; {\em mutatis
  mutandis,\/} the discussion which follows applies to this case as
well.

Concentrating on the first case, let
$$\Phi_i(y,z)\in \Cbb[[y]][z]$$
be a monic polynomial of degree $m_i$, defining an irreducible branch
of $\cC$ at $p$, not tangent to $y=0$. Then $\Phi_i$ splits
(uniquely) as a product of linear factors over the ring $\Cbb[[y^*]]$
of power series with {\em rational nonnegative\/} exponents:
$$\Phi_i(y,z)=\prod_{j=1}^{m_i} \left(z-f_{ij}(y)\right)\quad,$$
with each $f_{ij}(y)$ of the form
$$f(y)=\sum_{k\ge 0} \gamma_{\lambda_k} y^{\lambda_k}$$
with $\lambda_k\in \Qbb$, $1\le \lambda_0<\lambda_1<\dots$, and
$\gamma_{\lambda_k}\ne 0$. We call each such $z=f(y)$ a {\em formal
  branch\/} of $\cC$ at $p$. The branch is {\em tangent to
  $z=0$\/} if the dominating exponent $\lambda_0$ is~$>1$. The terms
$z-f_{ij}(y)$ in this decomposition are the Puiseux series for
$\cC$ at $p$.

We will need to determine $\lim_{t\to
  0}\cC\circ\alpha(t)$ as a union of `limits' of the individual
formal branches at $p$. The difficulty here resides in the fact that
we cannot perform an arbitrary change of variable in a power series
with fractional exponents. In the case in which we will need to do
this, however, $\alpha(t)$ will have the following special form:
$$\alpha(t)=\begin{pmatrix}
1 & 0 & 0 \\
t^a & t^b & 0 \\
r(t) & s(t)t^b & t^c
\end{pmatrix}$$
with $a<b\le c$ positive integers and $r(t)$, $s(t)$ polynomials
(satisfying certain restrictions, which are immaterial here).
The difficulty we mentioned may be circumvented by the following
{\em ad hoc\/} definition.

\begin{defin}\label{branchlimit} 
The {\em limit\/} of a formal branch $z=f(y)$, along a germ
$\alpha(t)$ as above, is defined by the dominant term in
$$(r(t)+s(t)t^by+t^cz)-f(t^a)-f'(t^a)t^by-f''(t^a)t^{2b}\frac{y^2}2
-\cdots$$
where $f'(y)=\sum \gamma_{\lambda_k}\lambda_k y^{\lambda_k-1}$ etc.
\end{defin}

By `dominant term' we mean the coefficient of the lowest power of
$t$ after cancellations. This coefficient is a polynomial in $y$ and $z$,
giving the limit of the branch according to our definition.

This definition behaves as expected: that is, the limit of $\cC$ is the 
union of the limits of its individual formal branches. This fact will be used
several times in the rest of the paper, and may be formalized as follows.

\begin{lemma}\label{limitlemma}
Let $\Phi(y,z)\in \Cbb[[y]][z]$ be a monic polynomial,
$$\Phi(y,z)=\prod_i \left(z-f_i(y)\right)$$
a decomposition over $\Cbb[[y^*]]$, and let $\alpha(t)$ have
the special form
above. Then the dominant term in $\Phi\circ\alpha(t)$ is the product
of the limits of the formal branches $z=f_i(y)$ along~$\alpha$, 
as in Definition~\ref{branchlimit}.\hfill\qed
\end{lemma}

Let $m$ be the multiplicity of $\cC$ at $p=(0,0)$.
For simplicity, we assume that no branches of $\cC$ are tangent
to the line $y=0$, leaving to the reader the necessary adjustments
in the presence of such branches.
We write the 
generator $F$ for the ideal of $\cC$ as a product of formal branches
$F =\prod_{i=1}^m (z-f_i(y))$. We will  focus on the formal branches
that are tangent to the line $z=0$, which may be written explicitly as
$$z=f(y)=\sum_{k\ge 0} \gamma_{\lambda_k} y^{\lambda_k}$$
with $\lambda_k\in \Qbb$, $1<\lambda_0<\lambda_1<\dots$, and
$\gamma_{\lambda_k}\ne 0$.

\smallskip
Now we begin the proof of Proposition~\ref{standardform}.

\subsection{Reduction to $q\ne 0$}
The first remark is that, under the assumptions that $q$, $r$, and $s$
do not all vanish, we may in fact assume that $q(t)$ is not identically zero.

\begin{lemma}\label{qnot0}
If $\alpha(t)$ is as in $(\dagger)$, and $q= 0$, then $\lim_{t\to 0}\mathcal
C\circ\alpha(t)$ is a rank-2 limit.
\end{lemma}

\begin{proof}
(Sketch.) Assume $q=0$, and study the action of $\alpha(t)$
on individual monomials $x^Ay^Bz^C$ in an equation for $\cC$:
$$m_{ABC}:=x^A y^B (r(t) x+s(t) t^b y+t^c z)^C t^{bB}\quad.$$
There are various possibilities for the vanishing of $r$ and $s$, but
the dominant terms in $m_{ABC}$ are always kernel stars, which are
rank-2 limits by Lemma~\ref{rank2lemma}.
\end{proof}

\subsection{Reduction to $b<c$}
By Lemma~\ref{qnot0} and the last part of Lemma~\ref{faber} we may replace 
$\alpha(t)$ with an equivalent germ
$$\begin{pmatrix}
1 & 0 & 0     \\
t^a & t^b & 0 \\
r(t) & s(t)t^b & t^c\end{pmatrix}$$
with $a<b\le c$, and $r(t)$, $s(t)$ polynomials of 
degree $<c$, $<(c-b)$ respectively and vanishing at $t=0$.

Next, we
have to show that if $\lim_{t\to 0}\cC\circ\alpha(t)$ is not a
rank-2 limit then $b<c$ and $r(t)$, $s(t)$ are as stated
in Proposition~\ref{standardform}.

\begin{lemma}\label{b=c}
Let $\alpha(t)$ be as above. If $b=c$, then $\lim_{t\to 0} 
\cC\circ\alpha(t)$ is a rank-2 limit.
\end{lemma}

\begin{proof}
Decompose $F(1:y:z)$ in $\Cbb[[y,z]]$: write $F(1:y:z)=G(y,z)\cdot H(y,z)$,
where $G(y,z)$ collects the branches that are {\em not\/} tangent to
$z=0$. If $b=c$, then necessarily~$s=0$:
$$\alpha(t)=\begin{pmatrix}
1 & 0 & 0 \\
t^a & t^b & 0 \\
r(t) & 0 & t^b
\end{pmatrix}\quad,$$
and further $a<v(r)$ (cf.~Lemma~\ref{faber}). 
The reader can verify that the limits of the branches collected 
in $G$ are supported on the kernel line $x=0$. The limit of each (formal)
branch collected in $H(y,z)$ may be computed as in 
Definition~\ref{branchlimit}, and is found to be given by a homogeneous 
equation in $x$ and $z$ only: that is, a $(0:1:0)$-star. 
It follows that the limit of $\cC$ is again a kernel star, hence a rank-2 
limit by Lemma~\ref{rank2lemma}.
\end{proof}

\subsection{End of the proof of Proposition~\ref{standardform}}\label{Eop}
By Lemma~\ref{b=c}, we may now assume 
that $\alpha$ is given by
$$\alpha(t)=\begin{pmatrix}
1 & 0 & 0\\
t^a & t^b & 0\\
r(t) & s(t)t^b & t^c
\end{pmatrix}$$
with the usual conditions on $r(t)$ and $s(t)$, and further $a<b<c$.

The limit of $\cC$ under $\alpha$ is analyzed  by studying limits of
formal branches. 

\begin{lemma}\label{otherbranches}
The limits of formal branches that
are not tangent to the line $z=0$ are necessarily $(0:0:1)$-stars.
Further, if $a<v(r)$, then the limit of such a branch is the kernel
line $x=0$.\hfill\qed
\end{lemma}

\begin{lemma}\label{tangentbranches}
The limit of a formal branch $z=f(y)$ tangent to the line $z=0$ 
is a $(0:0:1)$-star unless
\begin{itemize}
\item $r(t)\equiv f(t^a)\pmod{t^c}$;
\item $s(t)\equiv f'(t^a)\pmod{t^{c-b}}$.
\end{itemize}
\end{lemma}

\begin{proof}
The limit of the branch is given by the dominant terms in
$$r(t)+s(t)t^by+t^cz=f(t^a)+f'(t^a)t^by+\dots$$
If $r(t)\not\equiv f(t^a)\pmod{t^c}$, then the weight of the branch
is necessarily $<c$, so the ideal of the limit is generated by a
polynomial in $x$ and $y$, as needed. 
The same reasoning applies if $s(t)\not\equiv f'(t^a)\pmod{t^{c-b}}$.
\end{proof}

To verify the condition on $\frac ca$ stated in 
Proposition~\ref{standardform}, note that the limit of the formal branch 
$z=f(y)$ is now given by the dominant term in
$$r(t)+s(t) t^b y+t^c z=f(t^a)+f'(t^a)t^b y+\frac{f''(t^a)t^{2b}y^2}2
+\dots:$$
the dominant weight will be less than $c$ (causing the limit to be a
$(0:0:1)$-star) if $c>2b+v(f''(t^a))=2b+a(\lambda_0-2)$. The stated
condition follows at once, completing the proof of Proposition~\ref{standardform}.

\subsection{Characterization of type~V germs}\label{charaV}
In the following, we will replace $t$ by a root of $t$ in the germ obtained in
Proposition~\ref{standardform}, if necessary, in order to ensure that the 
exponents appearing in its expression are relatively prime integers; 
the resulting germ determines the same component of the~PNC. 

In order to complete the characterization of type~V germs
given in \S\ref{germlist}, we need to determine
the possible triples $a<b<c$ yielding germs contributing
components of the PNC. This determination is best performed
in terms of $B=\frac ba$ and $C=\frac ca$. Let
$$z=f(y)=\sum_{k\ge 0} \gamma_{\lambda_k} y^{\lambda_k}$$
with $\lambda_k\in \Qbb$, $1<\lambda_0<\lambda_1<\dots$, and
$\gamma_{\lambda_k}\ne 0$, be a formal branch tangent to $z=0$.
Every choice of such a branch and of a rational number $C=\frac ca>1$
determines a truncation 
$$f_{(C)}(y):=\sum_{\lambda_k<C} \gamma_{\lambda_k} y^{\lambda_k}\quad.$$

The choice of a rational number $B=\frac ba$ satisfying $1<B<C$
and $B\ge \frac{C-\lambda_0}2+1$
determines now a germ as prescribed by Proposition~\ref{standardform}:
$$\alpha(t)=\begin{pmatrix}
1 & 0 & 0     \\
t^a & t^b & 0 \\
\underline{f(t^a)} & \underline{f'(t^a) t^b} & t^c\end{pmatrix}$$
(choosing the smallest positive integer $a$ for which the entries of 
this matrix have integer exponents). Observe that the truncation
$\underline{f(t^a)}=f_{(C)}(t^a)$ is identically~0 if and only if
$C\le\lambda_0$. Also observe that $\underline{f'(t^a)t^b}$ is
determined by $f_{(C)}(t^a)$, as it equals the truncation to $t^c$ 
of~${(f_{(C)})}'(t^a)t^b$.

\begin{prop}\label{abc}
If $C\le\lambda_0$ or $B\ne \frac{C-\lambda_0}2+1$, then $\lim_{t\to
  0}\cC\circ\alpha(t)$ is a rank-2 limit.
\end{prop}

We deal with the different cases separately.

\begin{lemma}\label{Clelambda0}
If $C\le\lambda_0$, then $\lim_{t\to 0}\cC\circ\alpha(t)$ is a
$(0:1:0)$-star.
\end{lemma}

\begin{proof}
If $C=\frac ca\le\lambda_0$, then $f_{(C)}(y)=0$, so
$$\alpha(t)=\begin{pmatrix}
1 & 0 & 0     \\
t^a & t^b & 0 \\
0 & 0 & t^c\end{pmatrix}\quad.$$
The statement follows by computing the limit of individual formal branches,
using Definition~\ref{branchlimit}.
\end{proof}

By Lemma~\ref{rank2lemma}, the limits obtained in
Lemma~\ref{Clelambda0} are rank-2 limits, so the first part of
Proposition~\ref{abc} is proved. As for the second part,
the limit of a branch tangent to
$z=0$ depends on whether the branch truncates to $f_{(C)}(y)$ or not.
These cases are studied in the next two lemmas. Recall that, by our
choice, $B\ge \frac{C-\lambda_0}2+1$.

\begin{lemma}\label{nottrunc}
Assume $C>\lambda_0$,
and let
$z=g(y)$ be a formal branch tangent to $z=0$, such that $g_{(C)}(y)\ne
f_{(C)}(y)$. Then the limit of the branch is supported on a kernel line.
\end{lemma}

\begin{proof}
The limit of the branch is determined by the dominant terms in
$$\underline{f(t^a)}+\underline{f'(t^a)t^b}y+t^cz=g(t^a)+g'(t^a)
t^by+\dots .$$
As the truncations $g_{(C)}$ and $f_{(C)}$ do not agree, the dominant
term is independent of~$z$. Under our hypotheses on $B$ and $C$,
it is found to be independent of $y$ as well, as needed.
\end{proof}

\begin{lemma}\label{dominant}
Assume $C>\lambda_0$,
and let
$z=g(y)$ be a formal branch tangent to $z=0$, such that $g_{(C)}(y)=
f_{(C)}(y)$. Denote by $\gamma_C^{(g)}$ the coefficient of $y^C$ in
$g(y)$. 
\begin{itemize}
\item If $B> \frac{C-\lambda_0}2+1$, then the limit of the branch
  $z=g(y)$ by $\alpha(t)$ is the line
$$z=(C-B+1)\gamma_{C-B+1} y+\gamma_C^{(g)}\quad.$$

\item If $B= \frac{C-\lambda_0}2+1$, then the limit of the branch
  $z=g(y)$ by $\alpha(t)$ is the conic
$$z=\frac {\lambda_0(\lambda_0-1)}2\gamma_{\lambda_0}y^2+\frac
{\lambda_0+C}2\gamma_{\frac{\lambda_0+C}2}y+\gamma_C^{(g)}\quad.$$
\end{itemize}
\end{lemma}

\begin{proof}
Rewrite the expansion whose dominant terms give the limit of the
branch as:
$$t^c z=(g(t^a)-\underline{f(t^a)})+(g'(t^a)t^b-\underline{f'(t^a)t^b})y
+\frac{g''(t^a)}2 t^{2b} y^2+\dots$$
The dominant term has weight $c=Ca$ by our choices; if
$B> \frac{C-\lambda_0}2+1$ then the weight of the coefficient of $y^2$
exceeds $c$, so it does not survive the limiting process, and the limit
is a line. If $B= \frac{C-\lambda_0}2+1$, the term in $y^2$ is
dominant, and the limit is a conic. 
The explicit expressions given in the statement are obtained by
reading the coefficients of the dominant terms.
\end{proof}

We can now complete the proof of Proposition~\ref{abc}:
\begin{lemma}\label{completion}
Assume $C>\lambda_0$.
If $B>\frac{C-\lambda_0}2+1$, then the limit $\lim_{t\to 0} \cC
\circ\alpha(t)$ is a rank-2 limit.
\end{lemma}

\begin{proof}
We will show that the limit is necessarily a kernel star, which
gives the statement by Lemma~\ref{rank2lemma}.

As $B>1$, the coefficient $\gamma_{C-B+1}$ is determined by the
truncation $f_{(C)}$, and in particular it is the same for
all formal branches with that truncation. Since
$B>\frac{C-\lambda_0}2+1$, by Lemma~\ref{dominant} the branches
considered there contribute lines through the fixed point
$(0:1:(C-B+1)\gamma_{C-B+1})$. We are done if we check that all other
branches contribute a kernel line $x=0$: and this is implied by
Lemma~\ref{otherbranches} for branches that are not tangent to $z=0$
(note $a<v(r)$ for the germs we are considering), and by
Lemma~\ref{nottrunc} for formal branches $z=g(y)$ tangent to $z=0$ but
whose truncation $g_{(C)}$ does not agree with $f_{(C)}$.
\end{proof}

\subsection{Quadritangent conics}\label{quadritangent}
We are ready to complete the proof of Theorem~\ref{mainmain}, by determining
the limits of the last contributing germs. These have been reduced to
the form listed as type~V in \S\ref{germlist} (up to a coordinate change, 
and replacing $t$ by a root of $t$):
$$\alpha(t)=\begin{pmatrix}
1 & 0 & 0 \\
t^a & t^b & 0 \\
\underline{f(t^a)} & \underline{f'(t^a)t^b} & t^c
\end{pmatrix}$$
for some branch $z=f(y)=\gamma_{\lambda_0}y^{\lambda_0}+\dots$ of 
$\cC$ tangent to $z=0$ at $p=(0,0)$, and further satisfying
$C>\lambda_0$ and $B=\frac{C-\lambda_0}2+1$ for $B=\frac ba$, $C=\frac
ca$. Type~V components of the PNC will arise depending on the limit
$\lim_{t\to 0}\cC\circ\alpha(t)$, which we now determine.

\begin{lemma}\label{quadconics}
If $C>\lambda_0$ and $B=\frac{C-\lambda_0}2+1$, then the limit
$\lim_{t\to 0}\cC\circ\alpha(t)$ consists of a union of
quadritangent conics, with distinguished tangent equal to the kernel
line $x=0$, and of a multiple of the distinguished tangent line.
\end{lemma}

\begin{proof}
Both $\gamma_{\lambda_0}$ and $\gamma_{\frac{\lambda_0+C}2}$
are determined by the truncation $f_{(C)}$ (since $C>\lambda_0$);
hence the equations of the conics 
$$z=\frac {\lambda_0(\lambda_0-1)}2\gamma_{\lambda_0}y^2+\frac
{\lambda_0+C}2\gamma_{\frac{\lambda_0+C}2}y+\gamma_C$$
contributed (according to Lemma~\ref{dominant}) by different branches
with truncation $f_{(C)}$ may only differ in the coefficient $\gamma_C$.

It is immediately verified that all such conics are tangent to the
kernel line $x=0$, at the point $(0:0:1)$, and that any two distinct such
conics meet only at the point $(0:0:1)$; thus they are necessarily
quadritangent.

Finally, the branches that do not truncate to $f_{(C)}(y)$ must
contribute kernel lines, by Lemmas~\ref{otherbranches} and \ref{nottrunc}.
\end{proof}

The degenerate case in which only one conic arises corresponds to germs
not contributing components of the projective normal cone, by
dimension considerations. A component is present as
soon as there are two or more conics, that is, as soon as two branches
contribute distinct conics to the limit.

This leads to the description given in \S\ref{germlist}. We say that
a rational number $C$ is `characteristic' for $\cC$ (with
respect to $z=0$) if at least two formal branches of $\cC$
(tangent to $z=0$) have the same nonzero truncation, but different
coefficients for $y^C$.

\begin{prop}\label{typeV}
The set of characteristic rationals is finite.

The limit $\lim_{t\to 0}\cC \circ\alpha(t)$ obtained in
Lemma~\ref{quadconics} determines a component of the projective normal
cone precisely when $C$ is characteristic.
\end{prop}

\begin{proof}
If $C\gg 0$, then branches with the same truncation must in fact be
identical, hence they cannot differ at $y^C$, hence $C$ is not
characteristic. Since the set of exponents of any branch is discrete,
the first assertion follows.

The second assertion follows from Lemma~\ref{quadconics}: if
$C>\lambda_0$ and $B=\frac{C-\lambda_0}2+1$, then the limit is a union
of a multiple kernel line and conics with equation
$$z=\frac {\lambda_0(\lambda_0-1)}2\gamma_{\lambda_0}y^2+\frac
{\lambda_0+C}2\gamma_{\frac{\lambda_0+C}2}y+\gamma_C\quad:$$
these conics are different precisely when the coefficients $\gamma_C$
are different, and the statement follows.
\end{proof}

Proposition~\ref{typeV} leads to the procedure giving components of
type~V explained in \S\ref{details} (also cf.~\cite{MR2001h:14068},
\S2, Fact~5), concluding the proof of Theorem~\ref{mainmain}.

\section{Boundaries of orbits}\label{boundary}

We have now completed the set-theoretic description of the PNC determined
by an arbitrary plane curve~$\cC$. As we have argued in \S\ref{prelim}, this
yields in particular a description of the boundary of~$\OC$. In this section
we include a few remarks aimed at making this description more explicit.

If $\dim\OC=8$, then the boundary of $\OC$ consists of the image of
the union of the PNC and of the proper transform $R$ of the complement
of $\PGL(3)$ in $\Pbb^8$ (cf.~Remark~\ref{eluding}). Curves in the image
of $R$ are stars (Lemma~\ref{PNCtolimits}). Curves in the image of the 
components of the PNC 
belong to the orbit closures of the limits of the marker germs listed 
in \S\ref{germlist}. We have proved that this list is exhaustive; therefore, 
the boundary of a given curve $\cC$ may be determined (up to stars)
by identifying the marker germs for $\cC$, and taking the
union of the orbit closures of the (finitely many) corresponding limits.

This reduces the determination of the curves in the boundary of the orbit 
of a given curve to the determination for curves with {\em small\/} orbit
(i.e., of dimension~$\le 7$). We note that, for a curve $\cC$
with small orbit, some 
components of the PNC will in fact dominate $\OCbar$: indeed,
in this case $\cC$ has positive dimensional stabilizer in $\PGL(3)$;
the limit of a germ centered at a singular matrix and otherwise 
contained in the stabilizer is~$\cC$ itself. This germ can be chosen to
be equivalent to a marker germ, identifying a component of the PNC which 
dominates the orbit closure.

As mentioned in the introduction,
the boundary for a curve with small orbit may be determined by a direct
method. Indeed, for such a curve we have constructed 
in~\cite{MR2002d:14083} explicit sequences of 
blow-ups at nonsingular centers which resolve the indeterminacies of
the basic rational map, and hence dominate the corresponding 
graph~$\TPbb^8$. The boundary of the curve may be determined by 
studying the image in $\Pbb^N$ of the various exceptional divisors of
these blow-ups.

\begin{center}
\begin{figure}
\includegraphics[scale=.4]{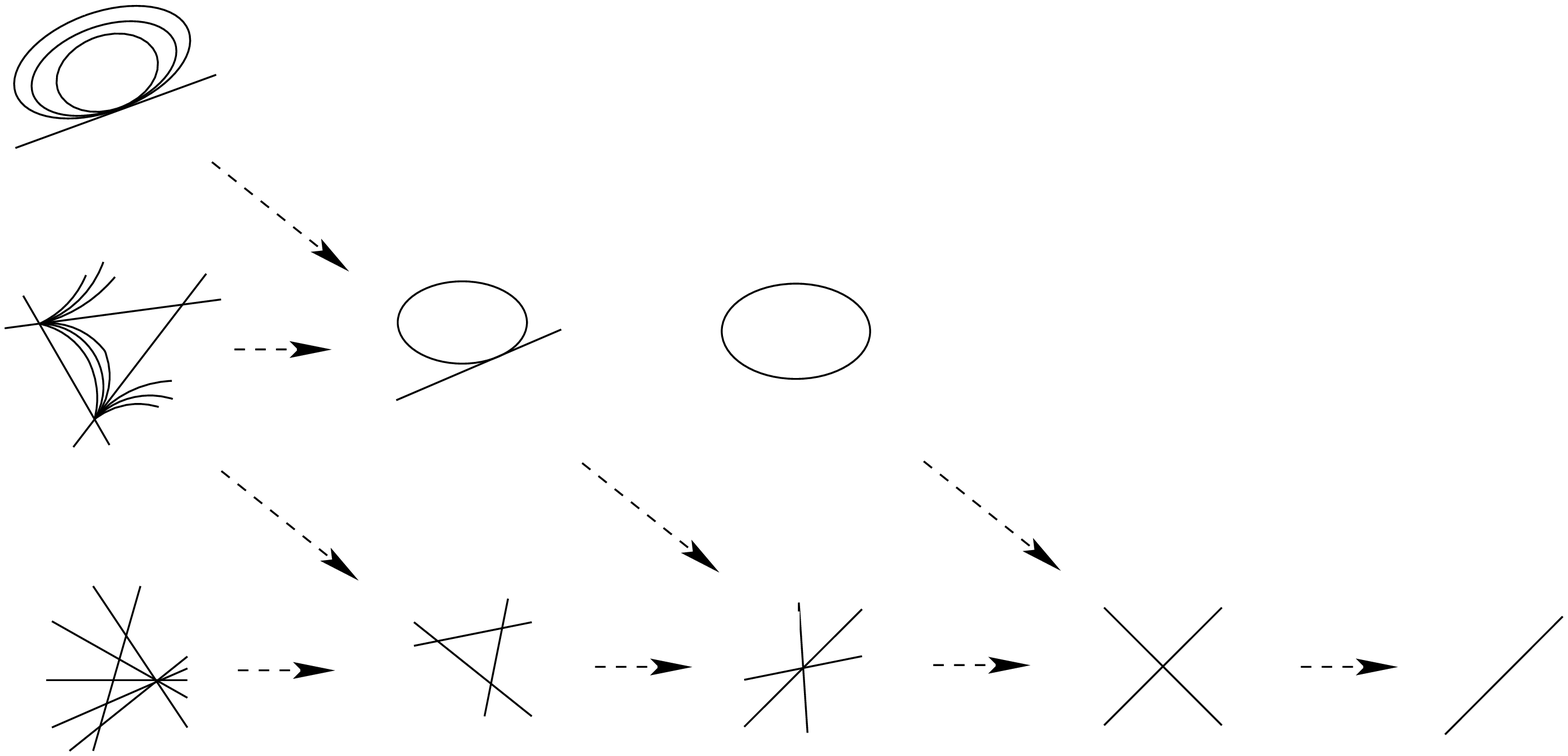}
\caption{}
\label{figure1}
\end{figure}
\end{center}

The result may be summarized by indicating which types of curves with
small orbits are in the boundary of a given curve with small orbit. 
Figure~\ref{figure1}
expresses part of this relation in terms of the representative
pictures for curves with small orbit shown in \S\ref{appendix}.
The five columns represent curves with orbits of dimension~7, 6, 5, 4, 2
respectively. Arrows indicate specialization: for example, 
the figure indicates that the boundary of the orbit of the union of a conic 
and a tangent line contains stars, but not single conics. Stars with more
than three lines are not displayed, to avoid cluttering the picture; 
the three kinds of curves displayed in the leftmost column all 
degenerate to such stars,
the only exceptions being the special cases of the second picture given
by the union of a conic and a transversal line, and by a single cuspidal cubic.

The situation illustrated here is precisely what one would expect
from naive considerations; it is confirmed by the study of
the blow-ups mentioned above.
Slightly more refined phenomena (for example, involving multiplicities
of the components) are not represented in this figure; in general,
they can be easily established by applying the results
of this paper or by analyzing the blow-ups of \cite{MR2002d:14083}. 

\begin{center}
\begin{figure}
\includegraphics[scale=.3]{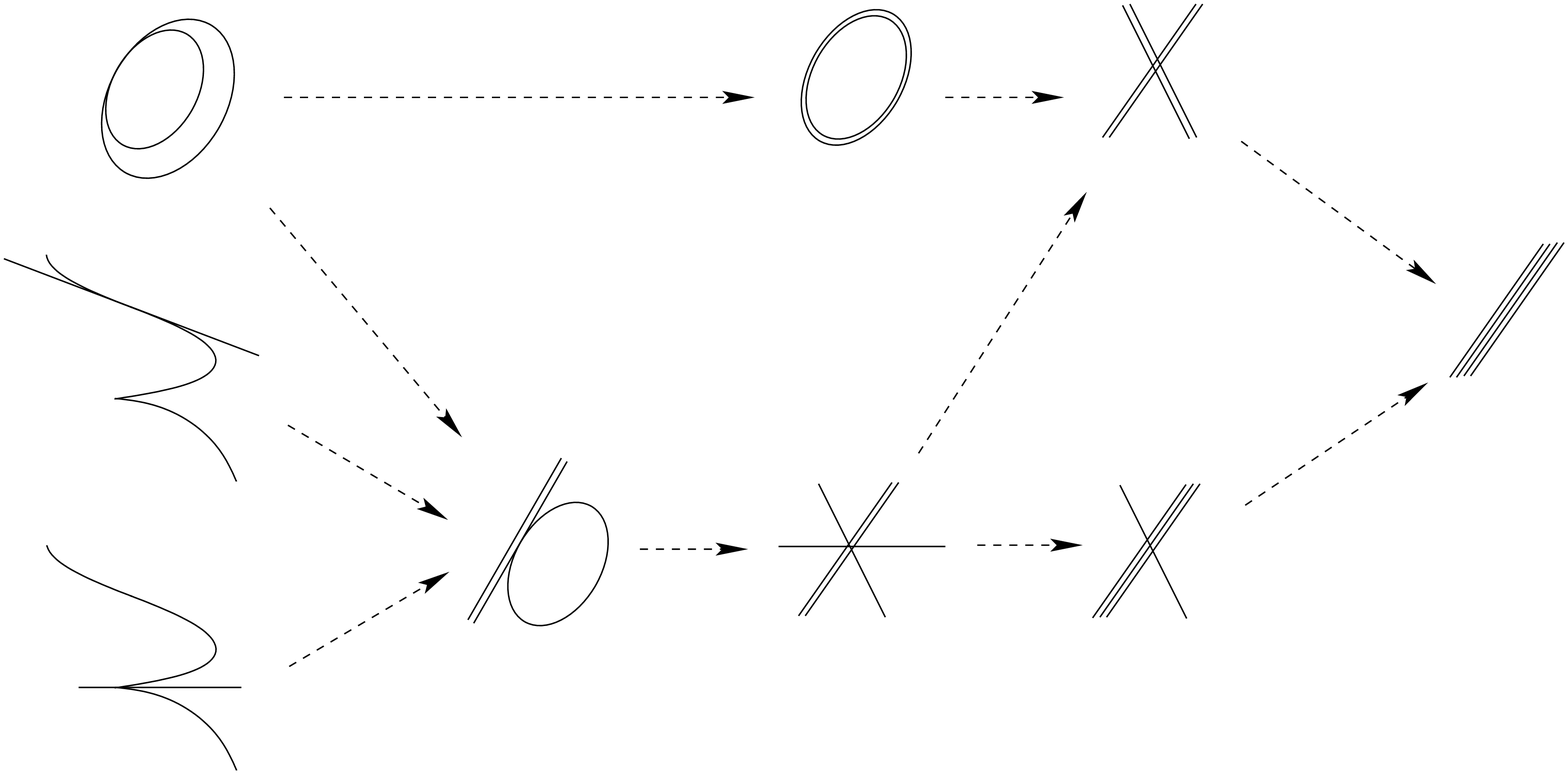}
\caption{}
\label{figure2}
\end{figure}
\end{center}

We close by pointing out one such phenomenon. In general, the union 
of a set of quadritangent conics
and a tangent line can specialize to the union of a conic and a
tangent line in two ways: (i) by type~II germs aimed at a general point of 
one of the conic components,
and (ii) by a suitable type~IV germ aimed at the tangency point.
The multiplicity of the conic in the limit is then the multiplicity
of the selected component in case (i), and the sum of the multiplicities 
of all conic components in case (ii). If the curve consists solely of quadritangent
conics, it degenerates to a multiple conic in case (ii). This possibility occurs
in the boundary of the orbit of the quartic curve from Example~\ref{extwo},
represented in Figure~\ref{figure2}.
We have omitted the set of stars of four distinct lines also in this
figure; in this case, it is a $6$-dimensional union of $5$-dimensional orbits.

\section{Appendix: curves with small linear orbits}\label{appendix}
For the convenience of the reader, we reproduce here the description of
plane curves with small linear orbits given in \cite{MR2002d:14084}. 
That reference contains a proof that this list is exhaustive, and details
on the stabilizer of each type of curve (as well as enumerative 
results for orbits of curves consisting of unions of lines, 
items (1)--(5) in the following list).

Let $\cC$ be a curve with small linear orbit. We list all possibilities for $\cC$,
together with the dimension of the orbit $\OC$ of $\cC$. The irreducible
components of the curves described here may appear with arbitrary
multiplicities.

\begin{enumerate}
\item $\cC$ consists of a single line; $\dim\OC=2$.
\item $\cC$ consists of 2 (distinct) lines; $\dim\OC=4$.
\item $\cC$ consists of 3 or more concurrent lines; $\dim\OC=5$. (We call
this configuration a {\it star\/}.)
\item $\cC$ is a triangle (consisting of 3 lines in general position);
$\dim\OC=6$.
\item $\cC$ consists of 3 or more concurrent lines, together with 1 other
(non-concur\-rent) line; $\dim\OC=7$. (We call this configuration a
{\it fan\/}.)
\item $\cC$ consists of a single conic; $\dim\OC=5$.
\item $\cC$ consists of a conic and a tangent line; $\dim\OC=6$.
\item $\cC$ consists of a conic and 2 (distinct) tangent lines;
$\dim\OC=7$.
\item $\cC$ consists of a conic and a transversal line and may contain either
one of the tangent lines at the 2 points of intersection or both of them;
$\dim\OC=7$.
\item $\cC$ consists of 2 or more bitangent conics (conics in the pencil
$y^2+\lambda x z$) and may contain the line $y$ through the two points
of intersection as well as the lines $x$ and/or $z$, tangent lines to the
conics at the points of intersection; again, $\dim\OC=7$.
\item $\cC$ consists of 1 or more (irreducible) curves from the pencil
$y^b+\lambda z^a x^{b-a}$, with $b\ge 3$, and may contain the lines $x$
and/or $y$ and/or $z$; $\dim\OC=7$.
\item $\cC$ contains 2 or more conics from a pencil through a conic and a
double tangent line; it may also contain that tangent line. In this case,
$\dim\OC=7$.
\end{enumerate}
The last case is the only one in which the maximal
connected subgroup of the stabilizer is the additive group $\Gbb_a$;
this fact was mentioned in \S\ref{germlist}.
The following picture represents schematically the curves described
above.
\begin{center}
\includegraphics[scale=.55]{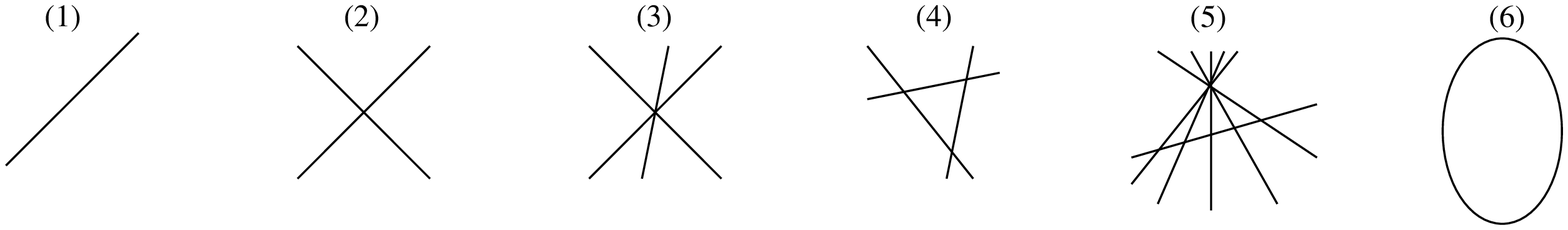}
\includegraphics[scale=.55]{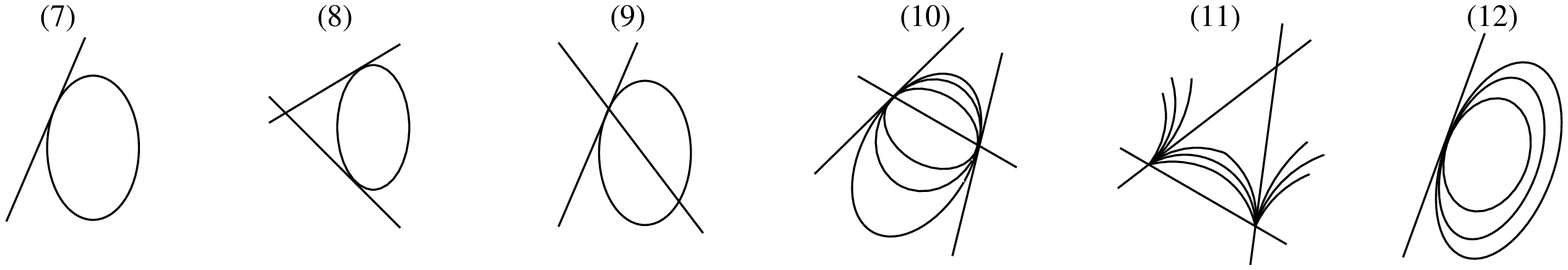}
\end{center}

\bibliographystyle{abbrv}
\bibliography{ghizzIbib}

\end{document}